# TESTING FOR COMMON ARRIVALS OF JUMPS FOR DISCRETELY OBSERVED MULTIDIMENSIONAL PROCESSES

By Jean Jacod and Viktor Todorov

*UPMC (Université Paris-6) and Northwestern University*

We consider a bivariate process $X_t = (X_t^1, X_t^2)$, which is observed on a finite time interval $[0, T]$ at discrete times $0, \Delta_n, 2\Delta_n, \ldots$. Assuming that its two components $X^1$ and $X^2$ have jumps on $[0, T]$, we derive tests to decide whether they have at least one jump occurring at the same time ("common jumps") or not ("disjoint jumps"). There are two different tests for the two possible null hypotheses (common jumps or disjoint jumps). Those tests have a prescribed asymptotic level, as the mesh $\Delta_n$ goes to 0. We show on some simulations that these tests perform reasonably well even in the finite sample case, and we also put them in use for some exchange rates data.

**1. Introduction.** It seems more and more apparent, as high-frequency data become available at a large scale, that many processes observed at discrete times, like stock prices or exchange rates, do have jumps. Now, finding models for discontinuous (continuous-time) processes that are compatible with data is a hard task, especially if one wants tractable models. This is even more difficult if one wants to model several processes at once.

Among models for multidimensional processes with correlated components, the easiest ones to tackle are those for which the various components do not jump together. Indeed, [16] assumes that jumps in individual stocks do not arrive together and can be diversified away when stocks are aggregated in a portfolio. But are such models for asset prices compatible with financial data? The empirical studies of [2] and [7], using high-frequency data, provide strong evidence for presence of jumps even on aggregate stock market level, which suggests that individual stocks contain a systematic jump component. It is clear that, if we want to formally study the systematic and idiosyncratic jumps in individual asset prices, we need formal tests for deciding whether the jumps in the different assets arrive together or not. Recently, [6] analyzed the relationship between the jumps in individual stocks and a portfolio of these stocks and similarly concluded for the need of formal tests about the possible common arrival of jumps in individual series. The main goal of this paper is to develop such tests in a general framework.

More specifically, we consider a $d$-dimensional process $X = (X^1, \ldots, X^d)$ that evolves according to a model, which we want to be as general as possible. We will take an Itô semimartingale, which essentially amounts to saying that it is driven by a Wiener process and a Poisson random measure (this allows in particular for "infinite activity" of the jumps). This semimartingale is observed at regularly spaced times $i\Delta_n$, where $i = 0, 1, \ldots$. We provide a testing procedure, based on the observation of the $X_{i\Delta_n}$'s up to some given terminal time $T$ (i.e., for $i = 0, 1, \ldots, [T/\Delta_n]$) to test the "null" hypothesis that two components, say $X^1$ and $X^2$, have no common jumps (meaning that they never jump together) on the time interval $[0, T]$, and also the null hypothesis that they do have common jumps.

An important feature of this paper is that we want these tests to be as independent of the underlying model as possible. An obvious second feature is that the problem is asymptotic; that is, the time lag $\Delta_n$ is "small" and, in fact, we study the asymptotic properties of the tests as $\Delta_n \to 0$, the horizon $T$ is kept fixed. An important third feature is that we test for common jumps or no, for the path of $t \mapsto X_t$ on $[0, T]$; some models allow for a positive probability of common jumps and simultaneously a positive probability for no common jump, and our tests try to give an answer for the observed path and not the model itself.

It might be useful to consider the case where the horizon $T = T_n$ also depends on $n$. There are two extreme cases: first, when $T_n$ converges to a limit $T > 0$, and all of what follows applies word for word in this case, although practically speaking this situation really amounts to considering $T_n = T$ as being constant; and second, when $T_n \to \infty$. Then, what follows does *not* apply, and we need much stronger assumptions on $X$, like ergodicity conditions, to derive any kind of results. On the other hand, the techniques are somewhat simpler and rather different, and we would be in a classical hypotheses testing situation instead of having "conditional" tests, as explained below. We do not consider this situation at all in this paper.

The tests exhibited here are based upon statistics involving suitable sums of functions of the increments of the process $X$ between successive observations. We will use these increments at two different scales, exactly as in [1], whose methods are generalized here. The way the tests are conducted is, however, different and in a sense more complicated than in that paper. Tests for deciding whether a given path has jumps or not on the interval $[0, T]$ have been already developed (see [1, 3, 4, 12, 13, 14, 15, 17] and [5] discuss some multivariate extensions of the test in [4]). Therefore, in this paper, we focus on the problem of testing whether there is at least one common jump time for the two components or none, supposing that there are jumps.

The paper is organized as follows. Sections 2 and 3 describe our setup and the test statistics we use. We provide a central limit theorem, or what plays the role of it for our proposed statistics, in Section 4, and use them to construct the actual tests in Section 5. We report the results of some Monte Carlo simulations in Section 6. In Section 7, we put our tests to use on actual data, namely the exchange rates between two pairs of currencies. Proofs are in Section 8.

## 2. Setting and assumptions.

Our problem in this paper is to determine whether any two components of a multidimensional process do jump at the same times. It thus amounts to solving this problem separately for each pair of components. In other words, this is a truly two-dimensional problem, and it is not a restriction to suppose that the underlying process $X$ is two-dimensional, with components denoted by $X^1$ and $X^2$.

As already mentioned, we do not want to make any specific model assumption on $X$, such as assuming some parametric family of models. We do need, however, a mild structural assumption that is satisfied in all continuous-time models with stochastic volatility used in finance, at least as long as one wants to rule out arbitrage opportunities.

Our structural assumption is that $X$ is an Itô semimartingale on some filtered space $(\Omega, \mathcal{F}, (\mathcal{F}_t)_{t \geq 0}, \mathbb{P})$, which means that it can be written as

$$
\begin{aligned}
(2.1) \quad X_t = X_0 &+ \int_0^t b_s \, ds + \int_0^t \sigma_s \, dW_s + \int_0^t \int \kappa \circ \delta(s,z)(\mu - \nu)(ds, dz) \\
&+ \int_0^t \int \kappa' \circ \delta(s,z)\mu(ds, dz),
\end{aligned}
$$

where $W$ and $\mu$ are a two-dimensional standard Wiener process and a Poisson random measure on $[0, \infty) \times E$, with $(E, \mathcal{E})$ an auxiliary measurable space, on the space $(\Omega, \mathcal{F}, (\mathcal{F}_t)_{t \geq 0}, \mathbb{P})$ and the predictable compensator (or intensity measure) of $\mu$ is $\nu(ds, dz) = ds \otimes \lambda(dz)$ for some given finite or $\sigma$-finite measure $\lambda$ on $(E, \mathcal{E})$. Above, $b$ is a two-dimensional adapted process, $\sigma$ is a $2 \times 2$-dimensional adapted process, and $\delta$ is a two-dimensional predictable function on $\Omega \times \mathbb{R}_+ \times E$. Moreover, $\kappa$ is a continuous truncation function on $\mathbb{R}^2$, that is a function from $\mathbb{R}^2$ into itself with compact support and $\kappa(x) = x$ on a neighborhood of 0, and $\kappa'(x) = x - \kappa(x)$.

Of course $b$, $\sigma$ and $\delta$ should be such that the integrals in (2.1) make sense (see, e.g., [9] for a precise definition of the last two integrals). However, we need a bit more than just the minimal integrability assumptions, and the precise hypotheses are stated in Assumption (H) below. Before this statement, we need some further notation. We write

$$(2.2) \qquad \Delta X_s = X_s - X_{s-}, \qquad \tau = \inf(t : \Delta X_t^1 \Delta X_t^2 \neq 0)$$

for the jumps of the $X$ process and the infimum $\tau$ of the joint jump times of the two components. Set also $\widetilde{\Gamma} = \{(\omega, t, x) : \delta^1(\omega, t, x)\delta^2(\omega, t, x) \neq 0\}$ and, for $i = 1, 2$,

$$(2.3) \qquad \delta_t'^{i}(\omega) = \begin{cases} \int (\kappa^i \circ \delta 1_{\widetilde{\Gamma}})(\omega, t, x)\lambda(dx), & \text{if the integral makes sense,} \\ +\infty, & \text{otherwise.} \end{cases}$$

ASSUMPTION (H).   (a) The paths $t \mapsto b_t(\omega)$ are locally bounded.

(b) The paths $t \mapsto \sigma_t(\omega)$ are all right-continuous with left limits.

(c) We have $\|\delta(\omega, t, x)\| \leq \Gamma_t(\omega)\gamma(x)$ identically, where $\Gamma$ is an adapted locally bounded process and $\gamma$ is a (nonrandom) nonnegative function satisfying $\int_E (\gamma(x)^2 \wedge 1)\lambda(dx) < \infty$.

(d) The paths $t \mapsto \delta_t'^{i}(\omega)$ for $i = 1, 2$ are locally bounded on the interval $[0, \tau(\omega))$.

(e) We have $\int_t^{t+r} \|\sigma_s\| \, ds > 0$ a.s. for all $t, r > 0$.

The nondegeneracy condition (e) says that, almost surely, the continuous martingale part of $X$ has no interval of constancy. It could be weakened, and, in any case, this condition is satisfied in all applications we have in mind. Apart from this nondegeneracy condition, which rules out "pure jump models" like the Variance Gamma or the NIG processes sometimes used in the financial literature, Assumption (H) accommodates virtually all models for stochastic volatility, including those with jumps, and allows for any kind of correlation or dependency between the volatility and asset price processes.

REMARK 2.1.   Condition (d) is implied by the others when $\int(\gamma(x) \wedge 1)\lambda(dx) < \infty$, which essentially amounts to saying that (and implies that) the jumps of $X$ are summable (i.e., $\sum_{s \leq t} \|\Delta X_s\| < \infty$ a.s.). Note that the summability of jumps in this sense implies at least that the processes $\delta'^{i}$ are locally bounded. Otherwise, it may appear as a strong assumption, because we can have $\delta_t'^{1}(\omega) = \delta_t'^{2}(\omega) = \infty$ for "most" $(\omega, t)$. However, if $A_t = \int_0^t \int 1_{\widetilde{\Gamma}}(s, x)\mu(ds, dx)$, then $A_\tau \leq 1$ by construction, and, by definition of the predictable compensator, we also have $\mathbb{E}(\int_0^\tau ds \int 1_{\widetilde{\Gamma}}(s, x)\lambda(dx)) = \mathbb{E}(A_\tau) \leq 1$. Hence, $\int_0^\tau |\delta_s'^{i}| \, ds < \infty$ a.s. Therefore, $|\delta_t'^{i}(\omega)| < \infty$ for $\mathbb{P}(d\omega) \otimes dt$-almost all $(\omega, t)$ such that $t \leq \tau(\omega)$. Hence, (d) is indeed a rather weak technical assumption, similar to saying that $b_t$ is locally bounded, instead of the "minimal" assumption saying that $\int_0^t \|b_s\| \, ds < \infty$ a.s.

*Comments on Assumption (H).* The key hypothesis is that $X$ is an Itô semimartingale; otherwise, everything falls apart. The nondegeneracy assumption (e) could perhaps be weakened to be $\int_0^T \|\sigma_s\| \, ds > 0$ a.s. for the final time $T$ only, but, without at least this weakened assumption, some of the forthcoming results are wrong. The conditions (a), (b) and (d), as seen in the previous remark, are rather weak, but they play an essential role in the proofs, and (b) also plays a crucial role in some of the statements. Finally, (c), which is not so weak, is also crucial for most proofs.

## 3. The test statistics.

3.1. *Preliminaries.* Recall that our process $X$ is observed over a given time interval $[0, T]$, at times $i\Delta_n$, for all $i = 0, 1, \ldots, [T/\Delta_n]$. We cannot, of course, do any better than if the process is observed "continuously" over $[0, T]$; that is, we can at best decide in which of the following three sets $\Omega_T^{(j)}$ (for "joint jumps"), $\Omega_T^{(d)}$ (for "disjoint jumps") or $\Omega_T^{(c)}$ (for "continuous"), the particular "observed" outcome $\omega$ lies:

$$(3.1) \quad \begin{cases} \Omega_T^{(j)} = \{\omega \colon \text{on } [0, T] \text{ the process } \Delta X_s^1 \Delta X_s^2 \text{ is not identically } 0\}, \\ \Omega_T^{(d)} = \{\omega \colon \text{on } [0, T] \text{ the processes } \Delta X_s^1 \text{ and } \Delta X_s^2 \text{ are not} \\ \qquad\qquad\quad \text{identically } 0, \text{ but the process } \Delta X_s^1 \Delta X_s^2 \text{ is}\}, \\ \Omega_T^{(c)} = \{\omega \colon \text{on } [0, T] \text{ at least one of } X^1 \text{ and } X^2 \text{ is continuous}\}. \end{cases}$$

That is, even under a "complete" observation of the path, we cannot decide whether the actual model allows for joint jumps or not, but only that the observed path has this property. Of course, if we decide that the observed path has joint jumps, then the model should allow for this; however, in the other case, the model can still allow for joint jumps.

These three sets are disjoint and form a partition of $\Omega$; however, we may very well have

$$(3.2) \qquad \mathbb{P}(\Omega_T^{(j)}) > 0, \qquad \mathbb{P}(\Omega_T^{(d)}) > 0, \qquad \mathbb{P}(\Omega_T^{(c)}) > 0,$$

at least in the case of finite activity jumps (e.g., when $\lambda$ is a finite measure). When both components have infinite activity, we have $\mathbb{P}(\Omega_T^{(c)}) = 0$, but the first two probabilities in (3.2) may still both be positive.

A comprehensive testing procedure should encompass all three kinds of outcomes. However, using the procedure established in [1] (or other methods), we can decide in principle whether we are in $\Omega_T^{(c)}$ or not. Here, we assume that this preliminary testing has been performed. If the conclusion is that we are in $\Omega_T^{(c)}$, then, of course, the procedure is ended. Otherwise, we have to decide between $\Omega_T^{(j)}$ and $\Omega_T^{(d)}$, which is the aim of this paper.

In a first case, we set the null hypothesis to be "joint jumps"; that is, we are in $\Omega_T^{(j)}$. We will take a critical (rejection) region $C_n^{(j)}$ at stage $n$, to be defined later, which should depend only on the observations $X_{i\Delta_n}$. Exactly as in [1], we do a kind of "conditional" test. Note that although $X$, and hence $\Omega_T^{(j)}$ as well, depend on the triple of coefficients $(b, \sigma, \delta)$ belonging to the set $\mathcal{H}$ of all coefficients satisfying Assumption (H), the observations $X_{i\Delta_n}$, and thus $C_n^{(j)}$, do *not* depend on $(b, \sigma, \delta)$ explicitly [the probability of $C_n^{(j)}$ does depend on this triple, though].

Then, with obvious notation, we take the following as our definition of the asymptotic size for a given triple of coefficients

$$(3.3) \qquad \alpha^{(j)} = \sup\Big(\limsup_n \mathbb{P}(C_n^{(j)} \mid A) : A \in \mathcal{F}, A \subset \Omega_T^{(j)}\Big).$$

Here, $\mathbb{P}(C_n^{(j)} \mid A)$ is the usual conditional probability, with respect to the set $A$, *with the convention that it vanishes if* $\mathbb{P}(A) = 0$. If $\mathbb{P}(\Omega_T^{(j)}) = 0$, then $\alpha^{(j)} = 0$, which is a natural convention, since, in this case, we want to reject the assumption whatever the outcome $\omega$ is. Note that $\alpha^{(j)}$ features some kind of "uniformity" over all subsets $A \subset \Omega_T^{(j)}$.

As for the asymptotic power, we define it as

$$(3.4) \qquad \beta^{(j)} = \inf\Big(\liminf_n \mathbb{P}(C_n^{(j)} \mid A) : A \in \mathcal{F}, A \subset \Omega_T^{(d)}\Big).$$

Again, this is a number. The asymptotic level and powers are defined here in a different way than in [1], where the level and the power were, respectively, $\alpha_0^{(j)} = \limsup_n \mathbb{P}(C_n^{(j)} \mid \Omega_T^{(j)}) \leq \alpha^{(j)}$ and $\beta_0^{(j)} = \liminf_n \mathbb{P}(C_n^{(j)} \mid \Omega_T^{(d)}) \geq \beta^{(j)}$. The results would be unchanged if we had taken $\alpha_0^{(j)}$ and $\beta_0^{(j)}$ as our definition.

In the second case, we set the null hypothesis to be "disjoint jumps"; that is, we are in $\Omega_T^{(d)}$. We take a critical region $C_n^{(d)}$ at stage $n$, again to be defined later, and the asymptotic size and power for the triple of coefficients $(b, \sigma, \delta)$ in $\mathcal{H}$ are

$$(3.5) \qquad \begin{cases} \alpha^{(d)} = \sup\Big(\limsup_n \mathbb{P}(C_n^{(d)} \mid A) : A \in \mathcal{F}, A \subset \Omega_T^{(d)}\Big), \\ \beta^{(d)} = \inf\Big(\liminf_n \mathbb{P}(C_n^{(d)} \mid A) : A \in \mathcal{F}, A \subset \Omega_T^{(j)}\Big). \end{cases}$$

3.2. *Construction of the critical regions.* We first need some notation. For any Borel function $f$ on $\mathbb{R}^2$, we write

$$(3.6) \qquad \Delta_i^n X = X_{i\Delta_n} - X_{(i-1)\Delta_n}, \qquad V(f, \Delta_n)_t = \sum_{i=1}^{[t/\Delta_n]} f(\Delta_i^n X).$$

Below, we also use $V(f, k\Delta_n)_t$ for $k$ an integer bigger than 1, meaning that we replace the stepsize $\Delta_n$ by $k\Delta_n$. That is, we have

$$(3.7) \qquad V(f, k\Delta_n)_t = \sum_{i=1}^{[t/k\Delta_n]} f(X_{ik\Delta_n} - X_{(i-1)k\Delta_n}).$$

Note that $V(f, \Delta_n)_t$, and also $V(f, k\Delta_n)_t$ for all $k \geq 2$, can be computed on the basis of the observations.

The following three functions will be of particular interest:

$$(3.8) \qquad f(x) = (x_1 x_2)^2, \qquad g_1(x) = (x_1)^4 \quad \text{and} \quad g_2(x) = (x_2)^4.$$

The critical regions $C_n^{(j)}$ and $C_n^{(d)}$ will be based upon the following two test statistics:

$$(3.9) \qquad \Phi_n^{(j)} = \frac{V(f, k\Delta_n)_T}{V(f, \Delta_n)_T} \quad \text{and} \quad \Phi_n^{(d)} = \frac{V(f, \Delta_n)_T}{\sqrt{V(g_1, \Delta_n)_T V(g_2, \Delta_n)_T}}.$$

Here, $k$ is an integer not less than 2 (typically $k = 2$ or $k = 3$), which is fixed throughout. Note that $\Phi_n^{(j)}$ depends on $k$, and both $\Phi_n^{(j)}$ and $\Phi_n^{(d)}$ depend on $T$.

The asymptotic behavior of these two statistics is crucial, and in order to give a description of it we need the notion of stable convergence in law, for which we refer, for example, to [9]. We also need some (cumbersome) further notation to describe the limits.

Recall that (H) is assumed. We denote by $(S_q)_{q \geq 1}$ a sequence of stopping times which exhausts the "jumps" of the Poisson measure $\mu$. Hence, for each $\omega$, we have $S_p(\omega) \neq S_q(\omega)$ if $p \neq q$, and that $\mu(\omega, \{t\} \times E) = 1$ if and only if $t = S_q(\omega)$ for some $q$. There are many ways of constructing those stopping times, but it turns out that what follows does not depend on the specific description of them. Next, we consider an auxiliary space $(\Omega', \mathcal{F}', \mathbb{P}')$ which supports a number of variables and processes:

- four sequences $(U_q)$, $(U_q')$, $(\overline{U}_q)$, $(\overline{U}_q')$ of two-dimensional $\mathcal{N}(0, I_2)$ variables;
- a sequence $(\kappa_q)$ of uniform variables on $[0, 1]$;
- a sequence $(L_q)$ of uniform variables on the finite set $\{0, 1, \ldots, k-1\}$, where $k \geq 2$ is some fixed integer;

and all these variables are mutually independent. Then, we put

$$(3.10) \qquad \widetilde{\Omega} = \Omega \times \Omega', \qquad \widetilde{\mathcal{F}} = \mathcal{F} \otimes \mathcal{F}' \quad \text{and} \quad \widetilde{\mathbb{P}} = \mathbb{P} \otimes \mathbb{P}'.$$

We extend the variables $X_t, b_t, \ldots$ defined on $\Omega$ and $U_n, \kappa_n, \ldots$ defined on $\Omega'$ to the product $\widetilde{\Omega}$ in the obvious way, without changing the notation. We write $\widetilde{\mathbb{E}}$ for the expectation with regard to $\widetilde{\mathbb{P}}$. Finally, we let $(\widetilde{\mathcal{F}}_t)$ be the smallest (right-continuous) filtration of $\widetilde{\mathcal{F}}$ containing the filtration $(\mathcal{F}_t)$ and such that $U_n$, $U'_n$, $\kappa_n$ and $L_n$ are $\widetilde{\mathcal{F}}_{S_n}$-measurable for all $n$. Obviously, $\mu$ is still a Poisson measure with compensator $\nu$, and $W$ is still a Wiener process on $(\widetilde{\Omega}, \widetilde{\mathcal{F}}, (\widetilde{\mathcal{F}}_t)_{t\geq 0}, \widetilde{\mathbb{P}})$. Finally, we define the two-dimensional variables

$$(3.11) \quad \begin{cases} R_q = \sqrt{\kappa_q}\sigma_{S_q-}U_q + \sqrt{1-\kappa_q}\sigma_{S_q}U'_q, \\ R'_q = \sqrt{L_q}\sigma_{S_q-}\overline{U}_q + \sqrt{k-1-L_q}\sigma_{S_q}\overline{U}'_q, \\ R''_q = R_q + R'_q. \end{cases}$$

Let us next define some auxiliary processes to be used sometimes in the forthcoming "laws of large numbers" and also later in the associated CLTs. As a rule, processes without "tilde" are defined on the original space $\Omega$, and those with "tilde" are on the extension $\widetilde{\Omega}$. Below, we write $c_t = \sigma_t\sigma_t^\star$ (the diffusion matrix of $X$). Then, we set

$$(3.12) \quad \begin{cases} B_t = \sum_{s\leq t}(\Delta X_s^1)^2(\Delta X_s^2)^2, & C_t = \int_0^t(c_s^{11}c_s^{22} + 2(c_s^{12})^2)\,ds, \\ B_t'^1 = \sum_{s\leq t}(\Delta X_s^1)^4, & B_t'^2 = \sum_{s\leq t}(\Delta X_s^2)^4, \end{cases}$$

$$(3.13) \quad \begin{cases} F_t = \frac{1}{2}\sum_{s\leq t}((\Delta X_s^1)^2(c_{s-}^{22} + c_s^{22}) + (\Delta X_s^2)^2(c_{s-}^{11} + c_s^{11})), \\ F'_t = 2\sum_{s\leq t}((\Delta X_s^1)^2(\Delta X_s^2)^4(c_{s-}^{11} + c_s^{11}) \\ \qquad\qquad + (\Delta X_s^1)^4(\Delta X_s^2)^2(c_{s-}^{22} + c_s^{22}) + 2(\Delta X_s^1\Delta X_s^2)^3(c_{s-}^{12} + c_s^{12})), \end{cases}$$

$$(3.14) \quad \begin{cases} \widetilde{D}_t = \sum_{q:\, S_q\leq t}((\Delta X_{S_q}^1 R_q^2)^2 + (\Delta X_{S_q}^2 R_q^1)^2), \\ \widetilde{D}''_t = \sum_{q:\, S_q\leq t}((\Delta X_{S_q}^1 R_q''^2)^2 + (\Delta X_{S_q}^2 R_q''^1)^2), \end{cases}$$

$$(3.15) \quad \begin{cases} \widetilde{G}_t = 2\sum_{q:\, S_q\leq t}((\Delta X_{S_q}^1)^2\Delta X_{S_q}^2 R_q'^2 + (\Delta X_{S_q}^2)^2\Delta X_{S_q}^1 R_q'^1), \\ \widetilde{G}'_t = 2\sum_{q:\, S_q\leq t}((\Delta X_{S_q}^1)^4(c_{S_q-}^{22} + c_{S_q}^{22})(R_q^2)^2 \\ \qquad\qquad + (\Delta X_{S_q}^2)^4(c_{S_q-}^{11} + c_{S_q}^{11})(R_q^1)^2). \end{cases}$$

The following theorem gives us the asymptotic behavior of our two test statistics, on the union $\Omega_T^{(j)} \cup \Omega_T^{(d)}$. As said before, we supposedly know that we are not in $\Omega_T^{(c)}$, so the behavior of the statistics on this set is of no importance for us. Recall that (H) is assumed throughout.

THEOREM 3.1. (a) *We have*

$$(3.16) \qquad \Phi_n^{(d)} \xrightarrow{\mathbb{P}} \begin{cases} B_T / \sqrt{B_T'^1 B_T'^2} > 0, & on \ \Omega_T^{(j)}, \\ 0, & on \ \Omega_T^{(d)}. \end{cases}$$

(b) *We have*

$$(3.17) \qquad \Phi_n^{(j)} \xrightarrow{\mathbb{P}} 1 \qquad on \ \Omega_T^{(j)}$$

*and $\Phi_n^{(j)}$ converges stably in law, in restriction to the set $\Omega_T^{(d)}$, to a variable, that is, a.s. different from 1 and given by*

$$(3.18) \qquad \widetilde{\Phi} = \frac{\widetilde{D}_T'' + kC_T}{\widetilde{D}_T + C_T}.$$

The last claim means that $\mathbb{E}(h(\Phi_n^{(j)})Y1_{\Omega_T^{(d)}}) \to \widetilde{\mathbb{E}}(h(\widetilde{\Phi})Y1_{\Omega_T^{(d)}})$ for all bounded $\mathcal{F}$-measurable variables $Y$ and all bounded continuous functions $h$ on $\mathbb{R}$. This is the definition of the stable convergence in law, in restriction to a subset of $\Omega$ (see [9] for more details on the stable convergence in law).

Of course, if either of the two sets $\Omega_T^{(j)}$ or $\Omega_T^{(d)}$ has a vanishing probability, the corresponding statement above is empty.

As a consequence, we are led to take critical regions of the form $C_n^{(j)} = \{|\Phi_n^{(j)} - 1| \geq \varepsilon_n\}$ or $C_n^{(d)} = \{\Phi_n^{(d)} \geq \varepsilon_n\}$ for suitable, and possibly random, sequences $\varepsilon_n$. However, to determine the level of such tests we need to go a bit further and give a central limit theorem associated with the convergences established in Theorem 3.1, at least in restriction to $\Omega_T^{(j)}$ for $\Phi_n^{(j)}$ and to $\Omega_T^{(d)}$ for $\Phi_n^{(d)}$.

## 4. Central limit theorems.

We have a genuine CLT for $\Phi_n^{(j)}$, on $\Omega_T^{(j)}$. We do not really have it for $\Phi_n^{(d)}$ on $\Omega_T^{(d)}$, but it is replaced by the stable convergence in law toward a positive random variable, similar to the convergence in (3.18).

The basic theorem, about nonstandardized statistics, goes as follows.

THEOREM 4.1. (a) *In restriction to the set $\Omega_T^{(j)}$, the sequence $\frac{\Phi_n^{(j)}-1}{\sqrt{\Delta_n}}$ converges stably in law to the variable $\widetilde{\Psi} = \widetilde{G}_T / B_T$, which, conditionally on $\mathcal{F}$, is centered with variance*

$$(4.1) \qquad \widetilde{\mathbb{E}}(\widetilde{\Psi}^2 \mid \mathcal{F}) = (k-1)F_T'/(B_T)^2$$

*and is even Gaussian conditionally on $\mathcal{F}$ if the processes $X$ and $\sigma$ have no common jumps.*

(b) *In restriction to the set* $\Omega_T^{(d)}$, *the sequences* $\frac{1}{\Delta_n}\Phi_n^{(d)}$ *converges stably in law to the positive variable* $\widetilde{\Phi}' = (\widetilde{D}_T + C_T)/\sqrt{B_T'^1 B_T'^2}$, *which, conditionally on* $\mathcal{F}$, *satisfies*

$$(4.2) \qquad \widetilde{\mathbb{E}}(\widetilde{\Phi}' \mid \mathcal{F}) = (F_T + C_T)/\sqrt{B_T'^1 B_T'^2}.$$

4.1. *Some consistent estimators.* To evaluate the level of tests based on the statistic $\Phi_n^{(j)}$ or $\Phi_n^{(d)}$, we need consistent estimators for the asymptotic mean or variance obtained in Theorem 4.1. That is, we need to estimate $F_T'$ and $B_T$, respectively, $F_T$, $C_T$, $B_T'^1$ and $B_T'^2$, on the set $\Omega_T^{(j)}$, respectively, $\Omega_T^{(d)}$.

For $B_T$, $B_T'^1$ and $B_T'^2$ a simple extension of [8], which is also used for the first part of (3.16) (see Section 8.3), gives us that

$$
(4.3) \qquad
\begin{aligned}
V(f, \Delta_n)_T &\xrightarrow{\mathbb{P}} B_T, \\
V(g_1, \Delta_n)_T &\xrightarrow{\mathbb{P}} B_T'^1, \\
V(g_2, \Delta_n)_T &\xrightarrow{\mathbb{P}} B_T'^2.
\end{aligned}
$$

For $C_T$ we can use multipower variations or truncated powers. This gives rise to the following two alternative estimators:

$$
\begin{aligned}
(4.4) \qquad \widehat{A}_T^n = \frac{\pi^2}{4\Delta_n} \sum_{i=1}^{[T/\Delta_n]-3} \Big( & |\Delta_i^n X^1 \Delta_{i+1}^n X^1 \Delta_{i+2}^n X^2 \Delta_{i+3}^n X^2| \\
& + \frac{1}{8} |\Delta_i^n(X^1+X^2)\Delta_{i+1}^n(X^1+X^2) \\
& \qquad \times \Delta_{i+2}^n(X^1+X^2)\Delta_{i+3}^n(X^1+X^2)| \\
& + \frac{1}{8} |\Delta_i^n(X^1-X^2)\Delta_{i+1}^n(X^1-X^2) \\
& \qquad \times \Delta_{i+2}^n(X^1-X^2)\Delta_{i+3}^n(X^1-X^2)| \\
& - \frac{1}{4} |\Delta_i^n(X^1+X^2)\Delta_{i+1}^n(X^1+X^2) \\
& \qquad \times \Delta_{i+2}^n(X^1-X^2)\Delta_{i+3}^n(X^1-X^2)| \Big),
\end{aligned}
$$

$$(4.5) \qquad \widehat{A}_T'^n = \frac{1}{\Delta_n} \sum_{i=1}^{[T/\Delta_n]} f(\Delta_i^n X) 1_{\{\|\Delta_i^n X\| \le \alpha \Delta_n^\varpi\}},$$

where, for the second one, we choose $\alpha > 0$ and $\varpi \in (0, \frac{1}{2})$ arbitrarily.

For $F_T$ and $F'_T$, things are more complicated, and we do as in [1] and take any sequence $k_n$ of integers satisfying

$$(4.6) \qquad k_n \to \infty, \qquad k_n \Delta_n \to 0$$

and then let $I_{n,-}(i) = \{i - k_n, i - k_n + 1, \ldots, i - 1\}$ if $i > k_n$ and $I_{n,+}(i) = \{i + 2, i + 3, \ldots, i + k_n + 1\}$ define two local windows in time of length $k_n \Delta_n$ just before and just after time $i\Delta_n$. Then, we set, for $i \geq 1 + k_n$ and $m, l$ equal to 1 or 2,

$$(4.7) \qquad \begin{cases} \widehat{c}(n,-)_i^{ml} = \dfrac{1}{k_n \Delta_n} \displaystyle\sum_{j \in I_{n,-}(i)} \Delta_j^n X^m \Delta_j^n X^l 1_{\{\|\Delta_j^n X\| \leq \alpha \Delta_n^\varpi\}}, \\[2ex] \widehat{c}(n,+)_i^{ml} = \dfrac{1}{k_n \Delta_n} \displaystyle\sum_{j \in I_{n,+}(i)} \Delta_j^n X^m \Delta_j^n X^l 1_{\{\|\Delta_j^n X\| \leq \alpha \Delta_n^\varpi\}}. \end{cases}$$

Those are "estimates" of the diffusion matrix $c_t$ on the left and on the right of time $i\Delta_n$, respectively. With this in mind, and with $I_n(i) = I_{n,-}(i) \cup I_{n,+}(i)$, the desired estimators are the following:

$$(4.8) \qquad \begin{aligned} \widehat{F}_t^n = \frac{1}{2k_n \Delta_n} \sum_{i=1+k_n}^{[t/\Delta_n]-k_n-1} \sum_{j \in I_n(i)} &((\Delta_i^n X^1)^2 (\Delta_j^n X^2)^2 \\ &+ (\Delta_i^n X^2)^2 (\Delta_j^n X^1)^2) \\ &\times 1_{\{\|\Delta_i^n X\| > \alpha \Delta_n^\varpi, \|\Delta_j^n X\| \leq \alpha \Delta_n^\varpi\}}, \end{aligned}$$

$$(4.9) \qquad \begin{aligned} \widehat{F}_t'^n = \frac{2}{k_n \Delta_n} \sum_{i=1+k_n}^{[t/\Delta_n]-k_n-1} \sum_{j \in I_n(i)} &(\Delta_i^n X^1)^2 (\Delta_i^n X^2)^2 \\ &\times (\Delta_i^n X^1 \Delta_j^n X^2 + \Delta_i^n X^2 \Delta_j^n X^1)^2 1_{\{\|\Delta_i^n X\| > \alpha \Delta_n^\varpi, \|\Delta_j^n X\| \leq \alpha \Delta_n^\varpi\}}. \end{aligned}$$

The following theorem establishes the behavior of these estimators.

THEOREM 4.2. *Let $\alpha > 0$ and $\varpi \in (0, 1/2)$.*
(a) *We have*

$$(4.10) \qquad \widehat{A}_T^n \xrightarrow{\mathbb{P}} C_T, \qquad \Delta_n \widehat{A}_T'^m \xrightarrow{\mathbb{P}} 0,$$

$$(4.11) \qquad \widehat{F}_T^n \xrightarrow{\mathbb{P}} F_T, \qquad \widehat{F}_T'^m \xrightarrow{\mathbb{P}} F'_T.$$

(b) *Moreover we have:*

$$(4.12) \qquad \widehat{A}_T'^m \xrightarrow{\mathbb{P}} C_T \qquad \text{on the set } \Omega_T^{(d)},$$

$$(4.13) \qquad \text{the sequence of variables } \left(\frac{1}{\Delta_n} \widehat{F}_T'^m 1_{\Omega_T^{(d)}}\right)_{n \geq 1} \text{ is tight.}$$

REMARK 4.3. One could prove that, in restriction to the set $\Omega_T^{(d)}$, the sequence of variables $\frac{1}{\Delta_n} \widehat{F}_T'^n$ converges stably in law to the variable $\widetilde{G}_T'$ defined by (3.15), but this fact is not used for our tests.

The above is not quite enough for deriving tests of given asymptotic size (these quantities give rise to tests with a size smaller, and often significantly smaller, than the prescribed level), except in case (a) of Theorem 4.1, when $X$ and $\sigma$ do not jump together. We need, in fact, a sort of "estimate" for the distribution of the variables $\widetilde{G}_T$ and $\widetilde{D}_T$ defined on the extended space, and conditionally on $\mathcal{F}$. For this, we first denote by $\widehat{\sigma}(n, \pm)_i$ an arbitrary (measurable) square-root of the matrix $\widehat{c}(n, \pm)_i$ in (4.7), and we define the two-dimensional variables

$$(4.14) \quad \begin{cases} R(n)_i = \sqrt{\kappa_i} \widehat{\sigma}(n, -)_i U_i + \sqrt{1 - \kappa_i} \widehat{\sigma}(n, +)_i U_i', \\ R'(n)_i = \sqrt{L_i} \widehat{\sigma}(n, -)_i \overline{U}_i + \sqrt{k - 1 - L_i} \widehat{\sigma}(n, +)_i \overline{U}_i' \end{cases}$$

[the variables $(\kappa_i, L_i, U_i, U_i', \overline{U}_i, \overline{U}_i')$ are the ones defined before (3.10)]. Finally, *on the extended space*, we define the following processes:

$$(4.15) \quad \widehat{D}_t^n = \sum_{i=1+k_n}^{[t/\Delta_n]-k_n-1} ((\Delta_i^n X^1 R(n)_i^2)^2 + (\Delta_i^n X^2 R(n)_i^1)^2) 1_{\{\|\Delta_i^n X\| > \alpha \Delta_n^{\varpi}\}},$$

$$(4.16) \quad \widehat{G}_t^n = 2 \sum_{i=1+k_n}^{[t/\Delta_n]-k_n-1} \Delta_i^n X^1 \Delta_i^n X^2 (\Delta_i^n X^1 R'(n)_i^2 + \Delta_i^n X^2 R'(n)_i^1)$$
$$\times 1_{\{\|\Delta_i^n X\| > \alpha \Delta_n^{\varpi}\}}.$$

THEOREM 4.4. *Assume that we have a sequence $Z_n$ of positive variables going in probability to some variable $Z > 0$, on the space $(\Omega, \mathcal{F}, \mathbb{P})$. Then*

$$(4.17) \quad \widetilde{\mathbb{P}}(|\widehat{G}_T^n| > Z_n \mid \mathcal{F}) \xrightarrow{\mathbb{P}} \widetilde{\mathbb{P}}(|\widetilde{G}_T| > Z \mid \mathcal{F}),$$

$$(4.18) \quad \widetilde{\mathbb{P}}(\widehat{D}_T^n > Z_n \mid \mathcal{F}) \xrightarrow{\mathbb{P}} \widetilde{\mathbb{P}}(\widetilde{D}_T > Z \mid \mathcal{F}).$$

4.2. *CLT for the standardized statistics.* Combining Theorems 4.1 and 4.2, and in view of the properties of the stable convergence in law, we immediately get the following (at this stage, we need no proof).

THEOREM 4.5. (a) *With*

$$(4.19) \quad \widehat{V}_n^{(j)} = \frac{\sqrt{\Delta_n(k-1)\widehat{F}_T'^n}}{V(f, \Delta_n)_T},$$

the variables $(\Phi_n^{(j)} - 1)/\widehat{V}_n^{(j)}$ converge stably in law, in restriction to the set $\Omega_T^{(j)}$, to a variable which, conditionally on $\mathcal{F}$, is centered with variance 1 and is even $\mathcal{N}(0,1)$ if the processes $X$ and $\sigma$ have no common jumps.

(b) With

$$
(4.20) \quad
\begin{aligned}
\widehat{V}_n^{(d)} &= \frac{\Delta_n(\widehat{F}_T^n + \widehat{A}_T^n)}{\sqrt{V(g_1, \Delta_n)_T V(g_2, \Delta_n)_T}}, \\
\widehat{V}_n'^{(d)} &= \frac{\Delta_n(\widehat{F}_T^n + \widehat{A}_T'^n)}{\sqrt{V(g_1, \Delta_n)_T V(g_2, \Delta_n)_T}},
\end{aligned}
$$

the variables $\Phi_n^{(d)}/\widehat{V}_n^{(d)}$ and $\Phi_n^{(d)}/\widehat{V}_n'^{(d)}$ converge stably in law, in restriction to the set $\Omega_T^{(d)}$, to a positive variable which, conditionally on $\mathcal{F}$, has expectation 1.

Another consequence of Theorems 4.1 and 4.4 is the following, which will be important for some of the tests later.

THEOREM 4.6. *Let $Z_n$ and $Z$ be as in Theorem 4.4, and $A \in \mathcal{F}$.*
(a) *If $A \subset \Omega_T^{(j)}$, we have*

$$
(4.21) \quad \mathbb{P}\left(A \cap \left\{\frac{|\Phi_n^{(j)} - 1| V(f, \Delta_n)_T}{\sqrt{\Delta_n}} > Z_n\right\}\right) \to \widetilde{\mathbb{P}}(A \cap \{|\widetilde{G}_T| > Z\}).
$$

(b) *If $A \subset \Omega_T^{(d)}$, and with either $\widehat{A}_n = \widehat{A}_T^n$ or $\widehat{A}_n = \widehat{A}_T'^n$, we have*

$$
(4.22) \quad
\begin{aligned}
&\mathbb{P}\left(A \cap \left\{\frac{\Phi_n^{(d)}\sqrt{V(g_1, \Delta_n)_T V(g_2, \Delta_n)_T}}{\Delta_n} > Z_n + \widehat{A}_n\right\}\right) \\
&\qquad \to \widetilde{\mathbb{P}}(A \cap \{\widetilde{D}_T > Z\}).
\end{aligned}
$$

**5. Testing for common jumps.** We now use the preceding results to construct actual tests, either for the null hypothesis that there are jumps but no common jumps for the two components of $X$, or for the null hypothesis that there are common jumps.

5.1. *When there are common jumps under the null hypothesis.* In a first case, we set the null hypothesis to be "common jumps," that is, we are in $\Omega_T^{(j)}$. For this, we use the test statistics $\Phi_n^{(j)}$ and, in view of (3.17) and (3.18), we associate the critical region of the form

$$
(5.1) \quad C_n^{(j)} = \{|\Phi_n^{(j)} - 1| \geq c_n^{(j)}\}
$$

for some sequence $c_n^{(j)} > 0$, possibly even a random sequence, but in that case depending only on the observations $X_{i\Delta_n}$.

As usual, we fix a level $\alpha \in (0,1)$ and wish to find $c_n^{(j)}$ so that (5.1) asymptotically achieves this level; that is, the level for which $\alpha^{(j)} \leq \alpha$ and, of course, $\alpha^{(j)} = \alpha$ if possible.

If we know that $X$ and $\sigma$ do not jump together, then (a) of Theorem 4.5 allows to achieve $\alpha^{(j)} = \alpha$, and we need the $\alpha$-absolute quantile of $\mathcal{N}(0,1)$; that is, the number $z_\alpha$ such that $\mathbb{P}(|U| \geq z_\alpha) = \alpha$ for a $\mathcal{N}(0,1)$ variable $U$. Otherwise we may rely on Bienaymé–Chebyshev inequality to construct a test for which $\alpha^{(j)} \leq \alpha$. Or, we can make use of (a) of Theorem 4.6 in the following way.

Recall that at stage $n$ we know the variables given by (4.7). Then, we can use a Monte-Carlo procedure to simulate $N_n$ copies of the variables $R(n)_i$ of (4.14) (i.e., we simulate $N_n$ copies of the variables $(\kappa_i, U_i, U_i')_{1 \leq i \leq [T/\Delta_n]}$, and use (4.14) and the same observed variables $\widehat{\sigma}(n\pm)_i$ always to compute the $R(n)_i$'s). Plugging these into (4.16), and again with the same observed increments $\Delta_i^n X$, we obtain $N_n$ copies $(\widehat{G}(j)_T^n : 1 \leq j \leq N_n)$ of the variable $\widehat{G}_T^n$. Then, we take the order statistics for the absolute values $|\widehat{G}_{n,1}| \geq |\widehat{G}_{n,2}| \geq \cdots \geq |\widehat{G}_{n,N_n}|$ for this family, and we set

$$(5.2) \qquad Z_n^{(j)}(\alpha) = |\widehat{G}_{n,[\alpha N_n]}|;$$

that is, the $\alpha$-absolute quantile of the empirical distribution of the family $(\widehat{G}(j)_T^n : 1 \leq j \leq N_n)$. With this notation, we construct three slightly different tests.

THEOREM 5.1.  (a) *Assume that the two processes $X$ and $\sigma$ do not jump together. If we set*

$$(5.3) \qquad c_n^{(j)} = z_\alpha \widehat{V}_n^{(j)},$$

*where $\widehat{V}_n^{(j)}$ is given by (4.19), then the asymptotic level of the critical region defined by (5.1) for testing the null hypothesis "common jumps" satisfies*

$$(5.4) \qquad \alpha^{(j)} \leq \alpha$$

*and, if further $\mathbb{P}(\Omega_T^{(j)}) > 0$, we have $\alpha^{(j)} = \alpha$ and even*

$$(5.5) \qquad A \subset \Omega_T^{(j)}, \qquad \mathbb{P}(A) > 0 \quad \Rightarrow \quad \mathbb{P}(C_n^{(j)} \,|\, A) \to \alpha.$$

(b) *If we set*

$$(5.6) \qquad c_n^{(j)} = \widehat{V}_n^{(j)}/\sqrt{\alpha},$$

*where $\widehat{V}_n^{(j)}$ is given by (4.19), then the asymptotic level of the critical region defined by (5.1) for testing the null hypothesis "common jumps" satisfies (5.4).*

(c) *Take a sequence $N_n \to \infty$. Define $Z_n^{(j)}(\alpha)$ by (5.2), and set*

$$(5.7) \qquad c_n^{(j)} = Z_n^{(j)}(\alpha) \frac{\sqrt{\Delta_n}}{V(f, \Delta_n)_T}.$$

*Then, the asymptotic level of the critical region defined by (5.1) for testing the null hypothesis "common jumps" satisfies (5.4). If further $\mathbb{P}(\Omega_T^{(j)}) > 0$ we have $\alpha^{(j)} = \alpha$ and even (5.5).*

Clearly, (a) is preferable if it can be used, and, otherwise, (c) is preferable. The choice of the sequence $N_n$ going to infinity is asymptotically arbitrary, but $n$ is given in practice, and $N_n$ should be big enough to have a good approximation of the "true" $\alpha$-quantile of the $\mathcal{F}$-conditional distribution of $\widehat{G}_T^n$. This of course depends on $\alpha$, and taking, for example, $N_n = 1000/\alpha$ seems to be a reasonable choice (if $\alpha = 0.05$, this means $N_n = 20{,}000$; this looks like a big number, but the simulation of our $N_n = 20{,}000$ copies take only a few seconds when the number of observations $[T/\Delta_n]$ is about 1000).

In the previous theorem, there is no statement about the power of the test, for a good reason. Indeed, using Remark 4.3, we can show that for any $A \subset \Omega_T^{(d)}$ with $\mathbb{P}(A) > 0$, and if such a set exists, then $\mathbb{P}(C_n^{(j)} \mid A)$ converges to a limit that is smaller than 1. However, a simple modification of the previous tests allows us to obtain the same results under the null hypothesis, *and a power equal to* 1 *under the alternative*. It goes as follows.

THEOREM 5.2. *Let $\alpha' > 0$ and $\varpi' \in (0, \frac{1}{2})$. Then, if we replace (5.3) and (5.6), respectively, by*

$$(5.8) \qquad \begin{aligned} c_n^{(j)} &= z_\alpha(\widehat{V}_n^{(j)} \wedge (\alpha' \Delta_n^{\varpi'})), \qquad \text{respectively,} \\ c_n^{(j)} &= \frac{1}{\sqrt{\alpha}}(\widehat{V}_n^{(j)} \wedge (\alpha' \Delta_n^{\varpi'})). \end{aligned}$$

*Then, the claims* (a) *and* (b) *of Theorem 5.1 hold; and, furthermore, in these cases, the asymptotic power is $\beta^{(j)} = 1$.*

Note that, in situations in which the continuous component is the dominant part of $X$, $\Phi_n^{(j)}$ is expected to be very "close" to $k$, while $\widehat{V}_n^{(j)}$ is expected to be "sufficiently" small [if there were no jumps in both series the sequence $(\frac{1}{\Delta_n}\widehat{F}_T^m)_{n \geq 1}$ converges to 0]. As a result, we expect (and this is later confirmed in the Monte Carlo) that the tests for common jumps, when the critical values are determined using Theorem 5.1, still have good power against alternatives in $\Omega_T^{(d)}$. Therefore, for practical purposes, the critical regions of Theorem 5.1 are probably sufficient.

5.2. *When there are no common jumps under the null hypothesis.* In a second case, we set the null hypothesis to be "no common jumps"; that is, we are in $\Omega_T^{(d)}$. We take the critical region to be

$$(5.9) \qquad C_n^{(d)} = \{\Phi_n^{(d)} \geq c_n^{(d)}\}$$

for some sequence $c_n^{(d)} > 0$.

Here, we have two ways for choosing $c_n^{(d)}$: first, we can use (b) of Theorem 4.5 and the Markov inequality; or second, we can use (b) of Theorem 4.6 and proceed as in the previous subsection, in which case we simulate $N_n$ copies of $\widehat{D}_T^n$, giving rise to the order statistics $\widehat{D}_{n,1} \geq \widehat{D}_{n,2} \geq \cdots \geq \widehat{D}_{n,N_n}$ of this family, and we set

$$(5.10) \qquad Z_n^{(d)}(\alpha) = \widehat{D}_{n,[\alpha/N_n]}.$$

THEOREM 5.3. (a) *If we set*

$$(5.11) \qquad c_n^{(d)} = \widehat{V}_n/\alpha,$$

*where $\widehat{V}_n$ is either $\widehat{V}_n^{(d)}$ or $\widehat{V}_n'^{(d)}$, as given by (4.20), then the asymptotic level and power of the critical region defined by (5.9) for testing the null hypothesis "no common jumps" satisfy*

$$(5.12) \qquad \alpha^{(d)} \leq \alpha, \qquad \beta^{(d)} = 1.$$

(b) *Take a sequence $N_n \to \infty$. Define $Z_n^{(d)}(\alpha)$ by (5.10), and put either $\widehat{A}_n = \widehat{A}_T^n$ or $\widehat{A}_n = \widehat{A}_T'^n$. Then, if*

$$(5.13) \qquad c_n^{(d)} = (Z_n^{(d)}(\alpha) + \widehat{A}_n) \frac{\Delta_n}{\sqrt{V(g_1,\Delta_n)_T V(g_2,\Delta_n)_T}},$$

*the asymptotic level and power of the critical region defined by (5.9) for testing the null hypothesis "no common jumps" satisfy (5.12). If further $\mathbb{P}(\Omega_T^{(d)}) > 0$, we have $\alpha^{(d)} = \alpha$ and even*

$$(5.14) \qquad A \subset \Omega_T^{(d)}, \qquad \mathbb{P}(A) > 0 \quad \Rightarrow \quad \mathbb{P}(C_n^{(d)} \mid A) \to \alpha.$$

Again, here (b) seems preferable to (a), and the simulation results given below *strongly* suggest that one should use (b).

REMARK 5.4. The test statistics above are insensitive to the scales used to measure $X^1$ and $X^2$. If we multiply each component $X^i$ by a constant $\lambda_i$, the test statistics are unchanged. However, this is not true of the standardized versions, because of the truncation $\alpha\Delta_n^{\varpi}$ that we use for the modulus $\|\Delta_i^n X\|$ of the increments. This presupposes that both components $X^1$ and $X^2$ have increments with roughly the same order of magnitude. If this is

not true, we should either first multiply the first component, say $X^1$, by a suitable constant in such a way that the averages of $|\Delta_i^n X^1|$ and of $|\Delta_i^n X^2|$ (or the averages after deleting, say, the 10% biggest increments) are close one to the other, or we can use two different levels of truncation; that is, replace the set $\{\|\Delta_i^n X\| \leq \alpha \Delta^n\}$ by $\{|\Delta_i^n X^1| \leq \alpha_1 \Delta_n^\varpi, |\Delta_i^n X^2| \leq \alpha_2 \Delta_n^\varpi\}$.

REMARK 5.5. In fact, the choice of the truncation level $\alpha \Delta_n^\varpi$, in order to put our tests in use, is a difficult one. Asymptotically, this choice does not matter, but in practice it does matter a lot. The idea is that $\alpha \Delta_n^\varpi$ should be slightly bigger than "most" of the increments when there is no jump, or no big jump; those increments being of order of magnitude $\|\sigma_t \Delta_i^n W\|$ with $(i-1)\Delta_n \leq t \leq i\Delta_n$. A good choice, supported by empirical evidence coming from simulation, seems to be $\varpi = 0.48$ or $\varpi = 0.49$, and $\alpha$ being of about 3 or 4 times the "average" value of $\|\sigma_t\|$. The latter is unknown, but, usually, one has a good idea of its order of magnitude.

**6. Simulation results.** In this section, we check the performance of our tests on simulated data. In the simulation study, we work with the simple model

$$dX_t^1 = X_t^1 \sigma_1 \, dW_t^1 + \alpha_1 \int_{\mathbb{R}} X_{t-}^1 x_1 \mu_1(dt, dx_1) + \alpha_3 \int_{\mathbb{R}} X_{t-}^1 x_3 \mu_3(dt, dx_3),$$

$$dX_t^2 = X_t^2 \sigma_2 \, dW_t^2 + \alpha_2 \int_{\mathbb{R}} X_{t-}^2 x_2 \mu_2(dt, dx_2) + \alpha_3 \int_{\mathbb{R}} X_{t-}^2 x_3 \mu_3(dt, dx_3),$$

where $\mathrm{cor}(W^1, W^2) = \rho$, the Poisson measures $\mu_1$, $\mu_2$ and $\mu_3$ are independent with compensators $\nu_i(dt, dx_i) = \lambda_i \frac{1_{(x_i \in [-h_i; -l_i] \cup [l_i; h_i])}}{2(h_i - l_i)} \, dt \, dx_i$ for $0 < l_i < h_i$ and $i = 1, 2, 3$. The initial values are $X_0^1 = X_0^2 = 1$. We did not make simulations with infinite activity jumps, but we consider different values of the jump intensities $\lambda_i$, the highest being 25. This is "almost" like infinite activity.

In the Monte Carlo study the observation length is one day (i.e., $T = 1$ day) consistent with the literature on testing for jumps in individual financial series. We simulate from the above-given process for a total of 5000 days. Since we are interested in the behavior of the tests on the sets $\Omega_T^{(j)}$ and $\Omega_T^{(d)}$, we discard days in the simulation on which there is no common and/or disjoint jump in the two series. On each day, we consider sampling $n = 100$, $n = 1600$ or $n = 25{,}600$ times, corresponding approximately to sampling every 5 minutes, 30 seconds or 1 second for a trading day of 6.5 hours or, equivalently, to sampling every 15 minutes, 1 minute or 3 seconds for a trading day of 24 hours. In each simulation, we compute the raw statistics $\Phi_n^{(j)}$ and $\Phi_n^{(d)}$ as well as their standardized versions, which are defined as

$$T_n^{(j)} = \frac{\Phi_n^{(j)} - 1}{\widehat{V}_n^{(j)}}, \qquad T_n^{(d)} = \frac{\Phi_n^{(d)}}{\widehat{V}_n'^{(d)}},$$

where we use the notation (4.19) and (4.20). $\Phi_n^{(j)}$ is computed with $k = 2$. For the calculation of $\widehat{V}_n^{(j)}$ and $\widehat{V}_n'^{(d)}$, we use a local window $k_n = 1/\sqrt{\Delta_n}$ and truncation level of $\alpha \Delta_n^\varpi = 0.03 \times \Delta_n^{0.49}$.

In Table 1, we report the parameter values for all cases considered. In all simulation scenarios $\sigma_1^2 = \sigma_2^2 = 8 \times 10^{-5}$, and, therefore, we do not report these values in the table. In all considered cases, the variance of the common and disjoint jumps is $2 \times 10^{-5}$. This leads to proportion of the jumps in the individual series total variation that is similar to one estimated from real financial data (see, e.g., [7]). Note that scenarios with higher on average number of (common or disjoint) jumps automatically imply that the jumps are of smaller size. The different parameter settings differ in the average number of jumps (resp., their size), whether jumps arrive together or not and in the correlation between the continuous components of the price.

On Figures 1 and 2, we plot the Monte Carlo distributions of the raw statistics $\Phi_n^{(d)}$ and $\Phi_n^{(j)}$ under the different scenarios. On Figures 3 and 4 we plot the rejection curves associated with the standardized tests $T_n^{(d)}$ and $T_n^{(j)}$ [i.e., the rejection rates of the tests for disjoint, resp., common jumps when the critical values of the tests are determined using Theorem 5.3, part (a), resp., Theorem 5.1, part (a)]. Finally, on Figure 5, we plot the rejection curves of the test for disjoint jumps when the critical values are computed using the simulation approach of Theorem 5.3, part (b). We also calculated the rejection rates for the test for common jumps when the critical values are determined using Theorem 5.1, part (c). These rates are very similar to the ones reported on Figure 4 and, therefore, we do not report them here. We summarize our findings from the Monte Carlo study as follows.

TABLE 1
*Parameter setting for the Monte Carlo*

| | | | | | | | | | | | | | |
|------|------|------------|-------|--------|------------|-------------|-------|--------|------------|------------|-------|--------|
| **Parameters** | | | | | | | | | | | | | |
| **Case** | $\rho$ | $\alpha_1$ | $\lambda_1$ | $l_1$ | $h_1$ | $\alpha_2$ | $\lambda_2$ | $l_2$ | $h_2$ | $\alpha_3$ | $\lambda_3$ | $l_3$ | $h_3$ |
| I-j | 0.0 | 0.00 | | | | 0.00 | | | | 0.01 | 1 | 0.05 | 0.7484 |
| II-j | 0.0 | 0.00 | | | | 0.00 | | | | 0.01 | 5 | 0.05 | 0.3187 |
| III-j | 0.0 | 0.00 | | | | 0.00 | | | | 0.01 | 25 | 0.05 | 0.1238 |
| I-m | 0.5 | 0.01 | 1 | 0.05 | 0.7484 | 0.01 | 1 | 0.05 | 0.7484 | 0.01 | 1 | 0.05 | 0.7484 |
| II-m | 0.5 | 0.01 | 5 | 0.05 | 0.3187 | 0.01 | 5 | 0.05 | 0.3187 | 0.01 | 5 | 0.05 | 0.3187 |
| III-m | 0.5 | 0.01 | 25 | 0.05 | 0.1238 | 0.01 | 25 | 0.05 | 0.1238 | 0.01 | 25 | 0.05 | 0.1238 |
| I-d0 | 0.0 | 0.01 | 1 | 0.05 | 0.7484 | 0.01 | 1 | 0.05 | 0.7484 | | | | |
| II-d0 | 0.0 | 0.01 | 5 | 0.05 | 0.3187 | 0.01 | 5 | 0.05 | 0.3187 | | | | |
| III-d0 | 0.0 | 0.01 | 25 | 0.05 | 0.1238 | 0.01 | 25 | 0.05 | 0.1238 | | | | |
| I-d1 | 1.0 | 0.01 | 1 | 0.05 | 0.7484 | 0.01 | 1 | 0.05 | 0.7484 | | | | |
| II-d1 | 1.0 | 0.01 | 5 | 0.05 | 0.3187 | 0.01 | 5 | 0.05 | 0.3187 | | | | |
| III-d1 | 1.0 | 0.01 | 25 | 0.05 | 0.1238 | 0.01 | 25 | 0.05 | 0.1238 | | | | |

- *Testing the null of common jumps.* Under the null hypothesis, $\Phi_n^{(j)}$ is concentrated around 1, with more dispersion from this value (and slight upward bias) for less frequent sampling and settings with higher number of (smaller) jumps. Under the alternative hypothesis of disjoint jumps, $\Phi_n^{(j)}$ is concentrated around 2 as expected from the result in Theorem 3.1, part (b). Comparing the case of zero correlation between the Brownian motions with that of perfect correlation, we see that $\Phi_n^{(j)}$ is more concentrated around 2 for the case of perfect correlation. Turning to the testing of the null of common jumps, we see that $T_n^{(j)}$ has the right size in all scenarios when the sampling frequency is $n = 25{,}600$. On the other hand, for the case of higher number of smaller jumps (i.e., cases III-j and III-m) and sampling frequency $n = 100$, the test over-rejects quite significantly. It is interesting to note that, for $n = 100$, the rejection curves under cases III-j, III-m, III-d0 and III-d1 (the dashed lines in the plots on the last row of Figure 4) look similar. This is clearly a finite sample problem (the solid lines, corresponding to sampling frequency $n = 25{,}600$ in the plots on the last row of Figure 4, behave as expected). The reason is that, for relatively low sampling frequency (e.g., $n = 100$), the small jumps are hard to disentangle from the Brownian moves. Finally, the last two columns of Figure 4 reveal that $T_n^{(j)}$ has very good power against all considered alternatives. This suggests that, for practical purposes, there is no need to truncate $\widehat{V}_n^{(j)}$ in the construction of the critical region, as in Theorem 5.2 (which was done to guarantee asymptotic power of the test for common jumps of 1).

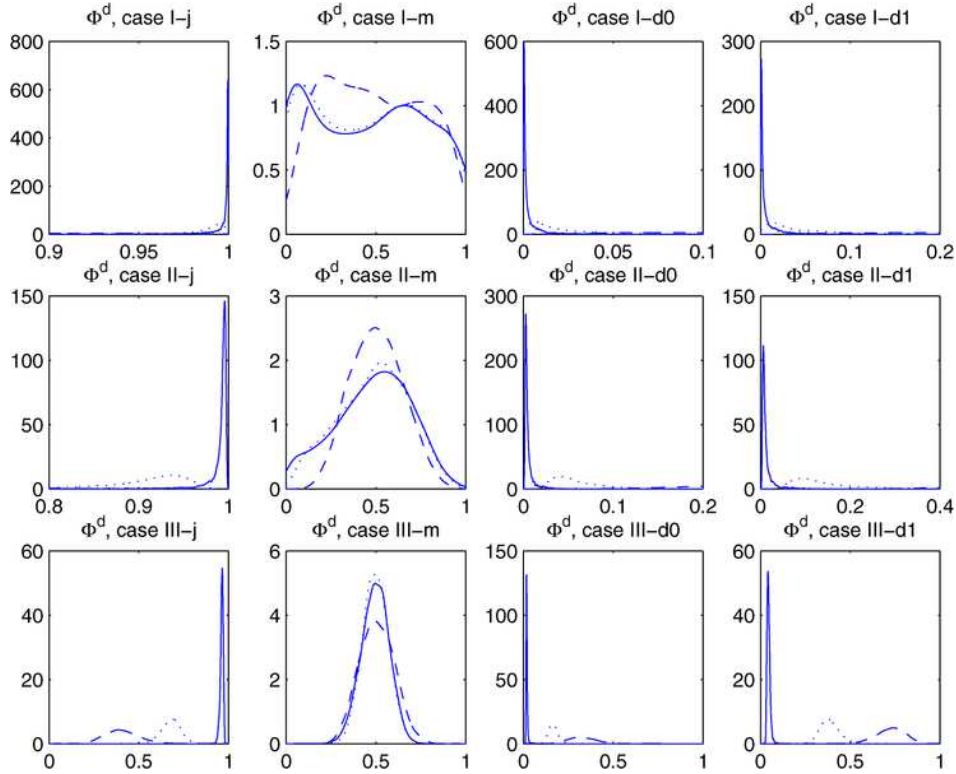

FIG. 1. *Kernel density estimate of $\Phi_n^{(d)}$ from the Monte carlo. The dashed line corresponds to sampling frequency of $n = 100$, the dotted line to sampling frequency of $n = 1600$ and the solid line to sampling frequency of $n = 25{,}600$.*

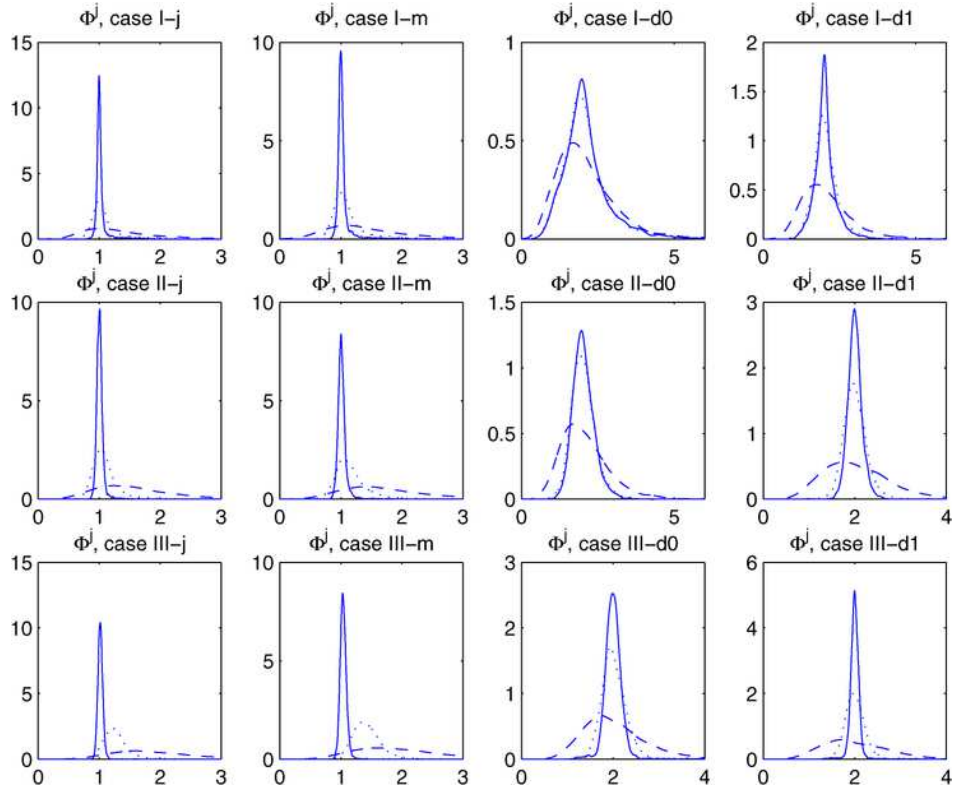

FIG. 2.   *Kernel density estimate of* $\Phi_n^{(j)}$ *from the Monte Carlo. The dashed line corresponds to sampling frequency of* $n = 100$, *the dotted line to sampling frequency of* $n = 1600$ *and the solid line to sampling frequency of* $n = 25{,}600$.

• *Testing the null of disjoint jumps.* Under the null hypothesis, consistently with the asymptotic results, $\Phi_n^{(d)}$ is concentrated around zero (see Figure 1, columns 3 and 4). Upward bias appears when the sampling frequency is low ($n = 100$), the number of jumps is higher (with smaller size) and the correlation between the Brownian motions in the prices is perfect.

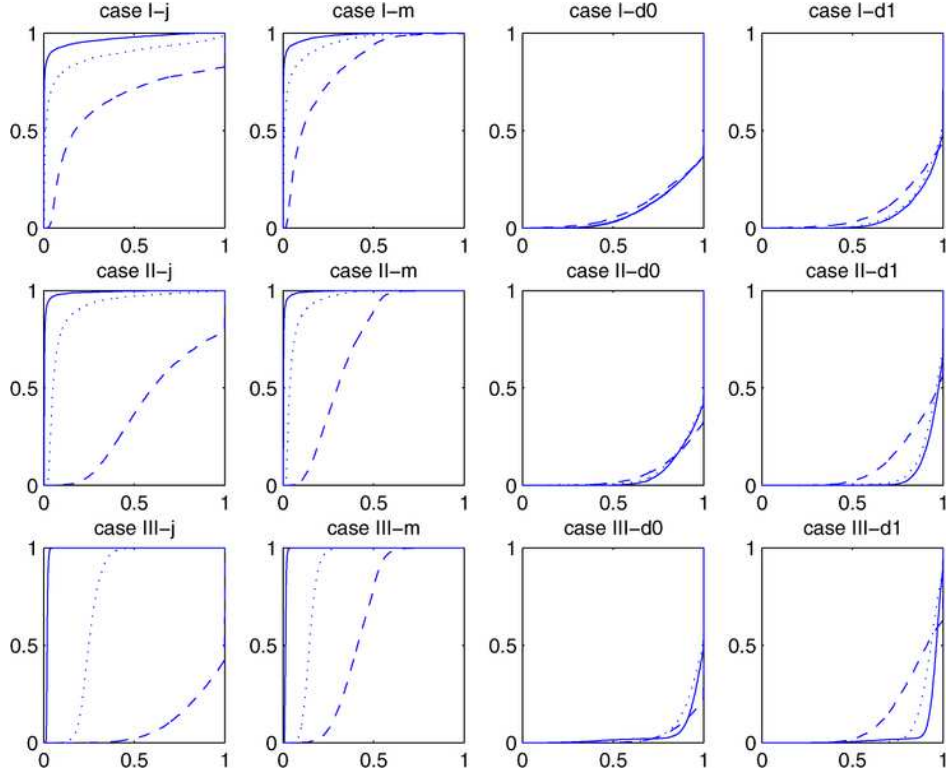

FIG. 3. *Size and power of the test for disjoint jumps with the critical values computed using Theorem 5.3, part* (a). *The x-axis shows the nominal level of the corresponding test, while the y-axis shows the percentage of rejection in the Monte Carlo. The dashed line corresponds to sampling frequency of $n = 100$, the dotted line to sampling frequency of $n = 1600$ and the solid line to sampling frequency of $n = 25{,}600$.*

Under the alternatives cases I-j, II-j and III-j, $\Phi_n^{(d)}$ takes values close to its asymptotic limit of 1. Under the alternatives cases I-m, II-m and III-m, $\Phi_n^{(d)}$ does not have a fixed nonrandom limit, because we have both common and disjoint jumps on the simulated trajectories. This is most clearly illustrated by the solid lines in the first and second plot of the second column of Figure 1. Note that the limiting value of $\Phi_n^{(d)}$ is almost surely different from 0 when we have both common and disjoint jumps. However, in a particular realization with relatively big disjoint jumps and small common jumps, the limiting value of $\Phi_n^{(d)}$ can get close to 0. Turning to the standardized test $T_n^{(d)}$, we see that using the Markov inequality in Theorem 5.3, part (a) leads to a significant underestimation of the size of the test. This holds true for all considered simulation scenarios. Therefore, it is recommendable to use the simulation approach in Theorem 5.3, part

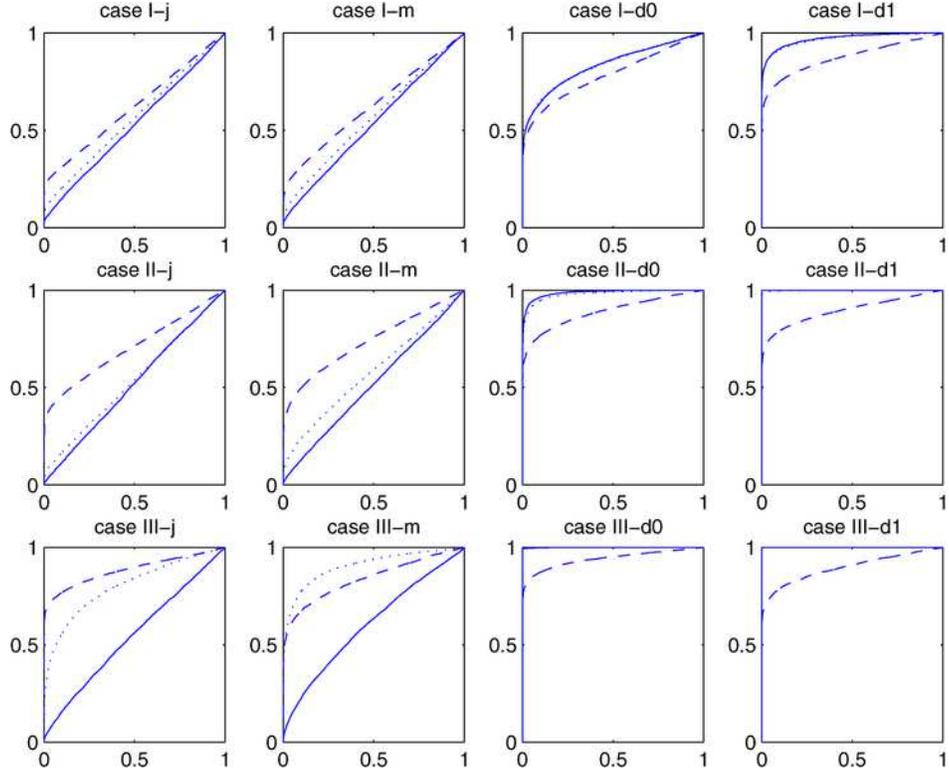

FIG. 4. *Size and power of the test for common jumps with critical values computed using Theorem 5.1, part* (a). *The x-axis shows the nominal level of the corresponding test, while the y-axis shows the percentage of rejection in the Monte Carlo. The dashed line corresponds to sampling frequency of $n = 100$, the dotted line to sampling frequency of $n = 1600$ and the solid line to sampling frequency of $n = 25{,}600$.*

(b) to determine the critical region of the test. As seen from Figure 5, when this is done, we do not have size distortions anymore. The only exception is case III-d1, where, even for $n = 25{,}600$, we have significant over-rejection. On the other hand, the first two columns of Figure 5 show that the test has very good power against the "common jumps" alternatives considered in the Monte Carlo. The only exception is for the lowest sampling frequency in cases III-j and III-m.

**7. Empirical application.** In the empirical part, we use high-frequency data from the foreign exchange spot markets for two exchange rates DM/$ and ¥/$. The data covers the period from December of 1986 through June of 1999, for a total of 3045 trading days. In each of the days, we sample every 5 minutes in the 24-hour trading day, and thus, in each of the trading days, we have 288 return observations. The DM/$ exchange rate data set

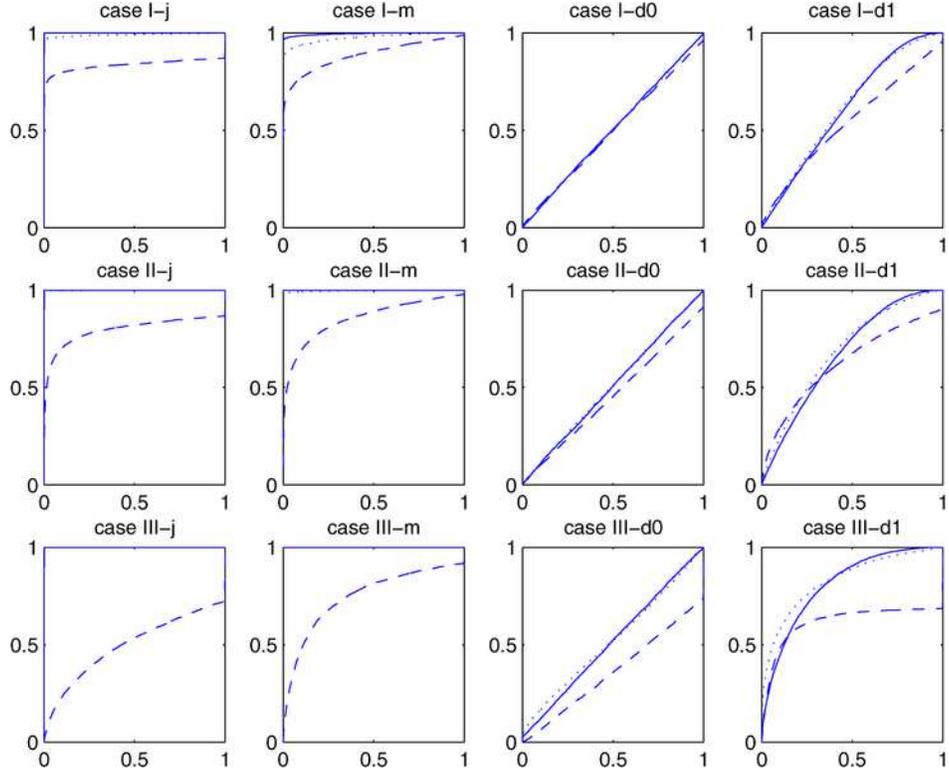

FIG. 5. *Size and power of the test for disjoint jumps with critical values computed via the simulation approach in Theorem 5.3, part* (b). *The x-axis shows the nominal level of the corresponding test, while the y-axis shows the percentage of rejection in the Monte Carlo. The dashed line corresponds to sampling frequency of $n = 100$, the dotted line to sampling frequency of $n = 1600$ and the solid line to sampling frequency of $n = 25{,}600$.*

has been used quite extensively in recent empirical studies for testing for presence of jumps (see, e.g., [2]). Here, we take this analysis one step further and test whether the jumps in the two studied exchange rate series arrive together.

First, we make several general comments on the empirical application of our tests. As mentioned in Section 3.1, on a first step, we need to remove the days in the sample (day is the time interval over which we apply our testing) on which at least one of the series does not exhibit jumps. This can be done with one of the many existing tests for presence of jumps in individual series. Once we select the days on which both series jump, we can perform our tests. We can apply them separately. For example, if we want to test the null of common jumps, we can use only our test statistic for common jumps $\Phi_n^{(j)}$. Alternatively, we can construct a rejection region for the common jump hypothesis by intersecting the rejection region of the test for common jumps (5.1) with the complement of the rejection region of the test for disjoint jumps (5.9). This is particularly attractive given the good power of both of our tests; that is, we can reduce the size with very little loss of power [as compared with the case when the testing is performed using only the critical region in (5.1)]. The same comments apply for the test for disjoint jumps. The particular choice of the critical values in constructing the rejection regions in the testing will depend on our tolerance toward Type I and Type II error.

We start the empirical analysis by identifying the days in the sample on which both series exhibit jumps. To be consistent with previous studies on the same data set (see [2]), we use the test based on the difference in logarithms of the realized variance and bi-power variation. The significance level of the test we chose is 1%. For this significance level, we found 288 days with jumps in the DM/$ exchange rate series and 291 days with jumps in the ¥/$ exchange rate series. Out of these days there are 40 days in which the tests indicate that both series contain jumps. In Table 2, we list these days, together with the raw statistics $\Phi_n^{(d)}$ and $\Phi_n^{(j)}$ and their $p$-values [which, following the conclusions from the Monte Carlo study, are computed using Theorem 5.3, part (b) and Theorem 5.1, part (c), resp.]. $\Phi_n^{(j)}$ is computed with $k = 2$. For the testing (i.e., calculation of the $p$-values) we set $T = 1$ (i.e., one day is our unit of measurement), and we used windows of size $k_n = \Delta_n^{-1/2} = 16$ and truncation level of $\alpha \Delta_n^\varpi = 3.0 \sqrt{\frac{BV(T)}{T}} \Delta_n^{0.49}$ applied to each individual series and the estimator $\hat{A}_T'^n$ of $C_T$. $BV(T)$, in the truncation level above, is the bi-power variation of the corresponding individual series and is a measure of $\int_0^T c_s^{jj} ds$ for $j = 1, 2$. We use it here to determine the magnitude of $c_t^{jj}$. Other alternative measures of $\int_0^T c_s^{jj} ds$ can also be used, and, to reduce the effect of measurement error, we can even use the whole sample period (and not just the day) to determine the level of $c_t^{jj}$.

Based on level of significance of 1%, we can separate the days in Table 2 in the following four categories.

- *Category* 1. The first category is of days in which the two tests find that there is common arrival of jumps. The total number of cojumping days is 22, which is significant, and the $p$-values associated with testing the null of disjoint jumps are very low (in most cases virtually 0) on these days.

- *Category* 2. The second category consists of days on which both tests indicate no common arrival of jumps. The number of these days is 5, which is small. However, note that we already found quite a significant number of days on which one of the series jumps while the other one does not. Importantly, the days in this category illustrate the possibility that, in spite of the fact that both series exhibit jumps during the day, they do so during different parts of the day. Note that if we were making our decision solely on the basis of individual tests for jumps, we would have misclassified the days in this category as days with common arrival of jumps (i.e., category 1). Thus, days in category 2 underline the importance of the tests developed here.

Table 2
*Empirical results for common jumps*

| Date | $\Phi_n^{(d)}$ | $\Phi_n^{(j)}$ | *p*-value $\Phi_n^{(d)}$ | *p*-value $\Phi_n^{(j)}$ |
|---|---|---|---|---|
| 09/11/1987 | 0.9938 | 1.0915 | 0.0000 | 0.4194 |
| 12/03/1987 | 0.6580 | 1.9831 | 0.0000 | 0.0342 |
| 12/10/1987 | 0.9933 | 1.1446 | 0.0000 | 0.2712 |
| 01/05/1988 | 0.5809 | 1.6876 | 0.0006 | 0.0276 |
| 01/15/1988 | 0.0040 | 1.6528 | 0.3663 | 0.5292 |
| 02/12/1988 | 0.9993 | 0.4100 | 0.0000 | 0.0038 |
| 05/17/1988 | 0.9658 | 1.0155 | 0.0000 | 0.8566 |
| 08/09/1988 | 0.5575 | 1.8825 | 0.0000 | 0.0404 |
| 09/14/1988 | 0.9984 | 0.7709 | 0.0000 | 0.2304 |
| 10/13/1988 | 0.9719 | 0.8011 | 0.0000 | 0.2792 |
| 10/26/1988 | 0.9731 | 1.3649 | 0.0000 | 0.1542 |
| 11/04/1988 | 0.9909 | 1.0527 | 0.0000 | 0.8476 |
| 05/17/1989 | 0.9860 | 0.6435 | 0.0000 | 0.1768 |
| 08/17/1989 | 0.9789 | 2.1938 | 0.0000 | 0.0002 |
| 09/27/1989 | 0.8255 | 1.1780 | 0.0000 | 0.6608 |
| 10/06/1989 | 0.9628 | 1.0647 | 0.0000 | 0.8320 |
| 10/17/1989 | 0.9732 | 1.4634 | 0.0000 | 0.1068 |
| 07/24/1991 | 0.8204 | 3.2959 | 0.0000 | 0.0002 |
| 08/02/1991 | 0.9753 | 1.2296 | 0.0000 | 0.3844 |
| 12/16/1991 | 0.2766 | 1.9990 | 0.0050 | 0.0002 |
| 01/10/1992 | 0.8595 | 0.6799 | 0.0000 | 0.3432 |
| 06/24/1992 | 0.9521 | 1.0435 | 0.0000 | 0.8692 |
| 08/24/1992 | 0.3306 | 2.0512 | 0.0018 | 0.0022 |
| 06/04/1993 | 0.9188 | 1.1350 | 0.0000 | 0.5880 |
| 09/16/1993 | 0.1866 | 1.2855 | 0.0222 | 0.6402 |
| 04/12/1994 | 0.2834 | 1.8069 | 0.0343 | 0.0088 |
| 06/17/1994 | 0.7766 | 2.6949 | 0.0000 | 0.0002 |
| 11/21/1994 | 0.1306 | 1.6013 | 0.3907 | 0.0834 |
| 03/17/1995 | 0.2787 | 2.7284 | 0.0267 | 0.0002 |
| 05/11/1995 | 0.6061 | 1.3020 | 0.0002 | 0.5118 |
| 11/13/1995 | 0.6948 | 2.2415 | 0.0000 | 0.0034 |
| 05/30/1996 | 0.5180 | 1.5381 | 0.0000 | 0.1440 |
| 06/27/1996 | 0.1544 | 0.7377 | 0.0010 | 0.4768 |
| 07/30/1997 | 0.1671 | 2.0925 | 0.7727 | 0.0004 |
| 03/30/1998 | 0.1203 | 2.4733 | 0.7621 | 0.0002 |
| 08/13/1998 | 0.1566 | 2.5072 | 0.2194 | 0.0006 |
| 10/05/1998 | 0.4035 | 1.4315 | 0.0164 | 0.1678 |
| 01/28/1999 | 0.1330 | 1.1790 | 0.0367 | 0.6524 |
| 03/01/1999 | 0.0498 | 1.9657 | 0.1661 | 0.0218 |
| 03/26/1999 | 0.2648 | 1.7011 | 0.0006 | 0.1178 |

- *Category* 3. The third category consists of days on which the null of both tests cannot be rejected—there are 6 days in this category. The possible explanations for such an outcome are at least two. First, it can be the case that on these days we have both common jumps and disjoint jumps with the magnitude of the common jumps far smaller as compared with the disjoint ones. The second possible explanation is that we have common jumps with very weak dependence. In both possible scenarios, the value of $\Phi_n^{(j)}$ is fairly close to 1, but for the above-mentioned reasons the value of $\Phi_n^{(d)}$ is close to zero. In general we will need more high-frequency observations for the $T_n^{(d)}$ test to gain power against such scenarios. Alternatively, we can perform the tests on different parts of the day.

- *Category* 4. The last category consists of days on which both tests reject their null hypothesis. The number of these days is 7. We notice that in these days the value of the $\Phi_n^{(d)}$ statistics is above 0.5 (i.e., it is relatively high), but the value of the $\Phi_n^{(j)}$ test is fairly close to 2.

Finally, if we use the less conservative significance level for determining the presence of jumps in the individual series of 5%, we find 113 days in which both series exhibit jumps. Further, if we use 5% significance level for testing for common and disjoint jumps we find that 55 of these days are in category 1, 10 in category 2, 11 in category 3 and 37 in category 4. Thus, our empirical study shows overall that the exchange series have a nontrivial number of days with common arrival of jumps, as well as days where jumps arrive at different times.

## 8. Proofs.

8.1. *Preliminaries.* We begin by showing that the processes $\widetilde{D}$, $\widetilde{D}''$, $\widetilde{G}$ and $\widetilde{G}'$ of (3.14) and (3.15) are actually well-defined and finite-valued and by stating some of their basic properties.

First, the process $\widetilde{D}$ is well-defined and increasing, but it might a priori take the value $+\infty$. However, by taking the $\mathcal{F}$-conditional expectation and using the properties of the variables $(\kappa_q, U_q, U_q')$, we get

$$\widetilde{\mathbb{E}}(\widetilde{D}_t \mid \mathcal{F}) = \sum_{q\,:\,S_q \leq t} \mathbb{E}'((\Delta X_{S_q}^1)^2 (R_q^2)^2 + (\Delta X_{S_q}^2)^2 (R_q^2)^2)$$

$$(8.1) \qquad = \tfrac{1}{2} \sum_{q\,:\,S_q \leq t} ((\Delta X_{S_q}^1)^2 (c_{S_q-}^{22} + c_{S_q}^{22}) + (\Delta X_{S_q}^2)^2 (c_{S_q-}^{11} + c_{S_q}^{11}))$$

$$= F_t.$$

See (3.13). Then, we deduce in particular that $\widetilde{D}$ is finite-valued, and the same argument shows that $\widetilde{D}''$ and $\widetilde{G}'$ are finite-valued as well. For $\widetilde{G}$, things are a bit more difficult, and we state the result in the form of a lemma (below, a process $U_t$, taking its values in the set of $2 \times 2$ nonnegative symmetric matrices, is said to be increasing for the strong order in this set if $U_t - U_s$ is a nonnegative matrix for all $s \leq t$).

LEMMA 8.1. *Let $\phi$ and $\psi$ be two real-valued functions on $\mathbb{R}^2$ with $\phi(x) = \mathrm{O}(\|x\|)$ and $\psi(x) = \mathrm{O}(\|x\|)$ as $x \to 0$.*

(a) *The process*

$$C(\phi, \psi)_t := \sum_{s \leq t} \phi(\Delta X_s)\psi(\Delta X_s)(c_{s-} + c_s) \tag{8.2}$$

*takes its values in the set of $2 \times 2$ nonnegative symmetric matrices, and it is increasing for the strong order in this set when $\psi = \phi$.*

(b) *The formulas (for $i = 1, 2$)*

$$Z^i(\phi)_t = \sum_{q : S_q \leq t} \phi(\Delta X_{S_q}) R_q^i, \qquad Z'^i(\phi)_t = \sum_{q : S_q \leq t} \phi(\Delta X_{S_q}) R_q'^i \tag{8.3}$$

*define two $\mathbb{R}^2$-valued processes $Z(\phi)$ and $Z'(\phi)$, and conditionally on $\mathcal{F}$ the eight-dimensional process $(Z(\phi), Z'(\phi), Z(\psi), Z'(\psi))$ is a square-integrable martingale with independent increments, zero mean and covariance given by*

$$\begin{cases} \widetilde{\mathbb{E}}(Z^i(\phi)_t Z^j(\psi)_t \mid \mathcal{F}) = \frac{1}{2}C(\phi, \psi)_t^{ij}, \\ \widetilde{\mathbb{E}}(Z'^i(\phi)_t Z^j(\psi)_t \mid \mathcal{F}) = 0, \\ \widetilde{\mathbb{E}}(Z'^i(\phi)_t Z'^j(\psi)_t \mid \mathcal{F}) = \dfrac{k-1}{2}C(\phi, \psi)_t^{ij}. \end{cases} \tag{8.4}$$

*Moreover, if $X$ and $c$ have no common jumps, the process $(Z(\phi), Z'(\phi), Z(\psi), Z'(\psi))$ is a Gaussian martingale, conditionally on $\mathcal{F}$.*

PROOF. This is proved exactly as Lemma 5.10 of [8]. The increasingness of $C(\phi, \phi)$ for the strong order comes from the fact that $c_s$ and $c_{s-}$ are nonnegative symmetric matrices. □

Then, with the notation (8.3), we obviously have $\widetilde{G} = Z'^1(f_1') + Z'^2(f_2')$ where $f_1'$ and $f_2'$ and the two first partial derivatives of the function $f$, and further with the notation (3.13),

$$\widetilde{\mathbb{E}}((\widetilde{G}_t)^2 \mid \mathcal{F}) = (k-1)F_t'. \tag{8.5}$$

Now, we state a strengthened version of Assumption (H).

ASSUMPTION (SH). We have Assumption (H), and $\|b_t\| + \|\sigma_t\| + \Gamma_t \leq K$ and also $\gamma(x) \leq K$ for some constant $K$. Then, up to multiplying $\gamma$ by a constant, we can even assume that $\|\delta(\omega, t, x)\| \leq \gamma(x)$.

If any of our limiting results holds under Assumption (SH), a *localization procedure* allows to get it under (H) only. This procedure, described in detail in [8], is omitted here. But in *all the remainder of the paper* we assume that Assumption (SH) holds.

We end up with some more notation. If $h$ is a function on $\mathbb{R}^2$, we have defined $V(h, \Delta_n)_t$ in (3.6), and when we want to emphasize the dependency on the process $X$, we write it $V(X; h, \Delta_n)_t$. We also use the following notation whenever the right-hand side below makes sense:

$$(8.6) \qquad V(X; h)_t = V(h)_t := \sum_{s \le t} h(\Delta X_s).$$

Next, we choose the functions $\psi$ on $\mathbb{R}_+$ and $\psi_a$ on $\mathbb{R}^2$, for $a > 0$, as follows:

$$(8.7) \quad \psi \text{ is decreasing, } C^\infty, \qquad 1_{[0,1]} \le \psi \le 1_{[0,2]}, \qquad \psi_a(x) = \psi(\|x\|/a).$$

Finally, we set

$$(8.8) \quad \gamma_0 = \sup_{z \in E} \gamma(z), \qquad A_\varepsilon = \{z : \gamma(z) \le \varepsilon\}, \qquad \Gamma'(\varepsilon) = \int_{A_\varepsilon} \gamma(z)^2 \lambda(dz).$$

8.2. *Estimates.* We derive some estimates for $X$, which follow from Assumption (SH) and will be used often in the sequel. Below, $K$ and $K_p$ are constants, changing from line to line, with $K_p$ depending on $p$. If $\varepsilon \in (0, \gamma_0]$, the process

$$(8.9) \qquad X(\varepsilon)_t = X_t - \int_0^t \int_{\{z : \gamma(z) > \varepsilon\}} \delta(s, z) \mu(ds, dz)$$

can be written as

$$X(\varepsilon)_t = X_0 + \int_0^t b(\varepsilon)_s \, ds + \int_0^t \sigma_s \, dW_s + \int_0^t \int_{A_\varepsilon} \delta(s, z)(\mu - \nu)(ds, dz),$$

where $b(\varepsilon)_t = b_t + \int (\delta(t, z) 1_{A_\varepsilon}(z) - \kappa \circ \delta(t, z)) \lambda(dz)$ is bounded by $K/\varepsilon$ for some $K$. We have $X(\varepsilon) = X_0 + X'(\varepsilon) + X''(\varepsilon)$, where

$$
(8.10) \quad
\begin{aligned}
X'(\varepsilon)_t &= \int_0^t b(\varepsilon)_s \, ds + \int_0^t \sigma_s \, dW_s, \\
X''(\varepsilon)_t &= \int_0^t \int_{A_\varepsilon} \delta(s, z)(\mu - \nu)(ds, dz).
\end{aligned}
$$

Note that $X(\gamma_0) = X$, so we will also use the notation

$$(8.11) \quad X' = X'(\gamma_0), \qquad X'' = X''(\gamma_0), \qquad \text{implying } X = X_0 + X' + X''.$$

First, $\|b(\varepsilon)\| \leq K/\varepsilon$ and the Davis–Burkholder–Gundy inequality yield for all $s, t \geq 0$ and $p \geq 1$:

$$
\begin{aligned}
(8.12) \quad & \mathbb{E}(\|X'(\varepsilon)_{s+t} - X'(\varepsilon)_s\|^p \mid \mathcal{F}_s) \\
& \leq K_p \left( t^{p/2} + \frac{t^p}{\varepsilon^p} \right),
\end{aligned}
$$

$$
\begin{aligned}
(8.13) \quad & \mathbb{E}(\|X'(\varepsilon)_{s+t} - X'(\varepsilon)_s - \sigma_s(W_{s+t} - W_s)\|^p \mid \mathcal{F}_s) \\
& \leq K_p \left( \frac{t^p}{\varepsilon^p} + \mathbb{E}\left( \left( \int_s^{s+t} \|\sigma_u - \sigma_s\|^2 \, du \right)^{p/2} \mid \mathcal{F}_s \right) \right).
\end{aligned}
$$

Next, if $p \geq 2$, the Davis–Burkholder–Gundy inequality and $\|\delta(t, z)\| \leq \gamma(z)$ yield that $\mathbb{E}(\|X''(\varepsilon)_{s+t} - X''(\varepsilon)_s\|^p \mid \mathcal{F}_s) \leq K_p \mathbb{E}(Z(\varepsilon)_t^{p/2})$ where $Z(\varepsilon)_t = \int_0^t \int_{A_\varepsilon} \gamma(z)^2 \mu(ds, dz)$ is a subordinator. Then, a well known result about Lévy processes (see, e.g., the proof of Lemma 5.1 of [10], with $H_t = 1$ and $a_t = t$), plus the obvious properties $\Gamma'(\varepsilon) \leq K$ and $\int_{A_\varepsilon} \gamma(z)^p \lambda(dz) \leq K\Gamma'(\varepsilon)$, give us $\mathbb{E}(Z(\varepsilon)_t^{p/2}) \leq Kt\Gamma'(\varepsilon)$ when $t \in [0, 1]$. Then, using Hölder inequality when $p < 2$, we get, for $s \geq 0$, $t \in [0, 1]$ and $p > 0$,

$$
(8.14) \quad \mathbb{E}(\|X''(\varepsilon)_{s+t} - X''(\varepsilon)_s\|^p \mid \mathcal{F}_s) \leq K_p (\Gamma'(\varepsilon) t)^{1 \wedge (p/2)}.
$$

Finally, using (45) of [1] and arguing componentwise, we obtain the existence of a an increasing function $\Gamma''$ on $\mathbb{R}_+$ with $\Gamma''(\varepsilon) \to 0$ as $\varepsilon \to 0$, such that for all $\eta > 0$ and $\theta \in (0, 1]$, and all $s \geq 0$ and $t \in [0, 1]$:

$$
(8.15) \quad \mathbb{E}(\|X''(\varepsilon)_{s+t} - X''(\varepsilon)_s\|^2 \wedge \eta^2 \mid \mathcal{F}_s) \leq Kt \left( \frac{\eta^2 + t}{\theta^2} + \Gamma''(\theta) \right).
$$

LEMMA 8.2. *Let $m, l \geq 2$ and let $j, k$ be two indices with values 1 or 2. Then, for all $\varepsilon \in (0, \gamma_0]$, we have*

$$
(8.16) \quad \mathbb{E}(|\Delta_i^n X'^j(\varepsilon)|^m |\Delta_i^n X''^k(\varepsilon)|^l \mid \mathcal{F}_{(i-1)\Delta_n}) \leq K\sqrt{\Gamma'(\varepsilon)} \Delta_n^2 \left( 1 + \frac{\Delta_n^{1/4}}{\varepsilon^m} \right).
$$

*If further the processes $X^1$ and $X^2$ have no common jump, then*

$$
(8.17) \quad \mathbb{E}(|\Delta_i^n X''^1(\varepsilon)|^m |\Delta_i^n X''^2(\varepsilon)|^l \mid \mathcal{F}_{(i-1)\Delta_n}) \leq K\sqrt{\Gamma'(\varepsilon)} \Delta_n^2.
$$

PROOF. (a) When $m > 2$, the estimate (8.16) is a simple consequence of (8.12), (8.14) and the Hölder inequality. When $m = 2$, an application of Itô's formula to (8.10) shows that

$$
(8.18) \quad |\Delta_i^n X'^j(\varepsilon)|^m |\Delta_i^n X''^k(\varepsilon)|^l = M(\varepsilon)_n + \int_{(i-1)\Delta_n}^{i\Delta_n} h(n, \varepsilon)_s \, ds,
$$

where $\mathbb{E}(M(\varepsilon)_n \mid \mathcal{F}_{(i-1)\Delta_n}) = 0$ and, with the notation $U_t = X'^j(\varepsilon)_t - X'^j(\varepsilon)_{(i-1)\Delta_n}$ and $V_t = X''^k(\varepsilon)_t - X''^k(\varepsilon)_{(i-1)\Delta_n}$ for $t \geq (i-1)\Delta_n$,

$$h(n,\varepsilon)_s = 2U_s|V_s|^l b(\varepsilon)_s^j + |V_s|^l c_s^{jj}$$
$$+ |U_s|^2 \int_{A_\varepsilon} (|V_{s-} + \delta^k(s,z)|^l - |V_{s-}|^l - l\{V_{s-}\}^{l-1}\delta^k(s,z))\lambda(dz)$$

(above, $\{v\}^r = |v|^r \mathrm{sign}(v)$ for any $v \in \mathbb{R}$). The integrand above is smaller than $K(1 + |V_{s-}|^{l-2})\gamma(z)^2$. Then, since $c_t$ is bounded and $|b(\varepsilon)_t^j| \leq K/\varepsilon$, we deduce from (8.12), (8.14), Hölder inequality and $\Gamma'(\varepsilon) \leq K$ that

$$\mathbb{E}(|h(n,\varepsilon)_{(i-1)\Delta_n+s}| \mid \mathcal{F}_{(i-1)\Delta_n}) \leq Ks\sqrt{\Gamma'(\varepsilon)}\left(1 + \frac{s^{1/4}}{\varepsilon^2}\right),$$

when $s \leq 1$, and (8.16) follows.

(b) For $t \geq (i-1)\Delta_n$, we write $U_t = X''^1(\varepsilon)_t - X''^1(\varepsilon)_{(i-1)\Delta_n}$ and $V_t = X''^2(\varepsilon)_t - X''^2(\varepsilon)_{(i-1)\Delta_n}$. Itô's formula yields that $|\Delta_i^n X''^1(\varepsilon)|^m|\Delta_i^n X''^2(\varepsilon)|^l$ equals the right-hand side of (8.18), where $M(\varepsilon)_n$ still has a vanishing conditional expectation, and $h(n,\varepsilon)_s$ has the form $h(n,\varepsilon)_s = \int_{A_\varepsilon} \alpha_{n,\varepsilon}(s,z)\lambda(dz)$, where

$$\alpha_{n,\varepsilon}(s,z) = |U_{s-} + \delta^1(s,z)|^m|V_{s-} + \delta^2(s,z)|^l - |U_{s-}|^m|V_{s-}|^l$$
$$- m\{U_{s-}\}^{m-1}|V_{s-}|^l\delta^1(s,z) - l|U_{s-}|^m\{V_{s-}\}^{l-1}\delta^2(s,z).$$

Now, if $X^1$ and $X^2$ never jump together, the product $\delta^1\delta^2$ vanishes $\mathbb{P}(d\omega) \otimes ds \otimes \lambda(dz)$ almost everywhere. Hence, $\alpha_{n,\varepsilon}$ is almost everywhere equal to

$$\alpha'_{n,\varepsilon}(s,z) = |V_{s-}|^l(|U_{s-} + \delta^1(s,z)|^m - |U_{s-}|^m - m\{U_{s-}\}^{m-1}\delta^1(s,z))$$
$$+ |U_{s-}|^m(|V_{s-} + \delta^2(s,z)|^l - |V_{s-}|^l - l\{V_{s-}\}^{l-1}\delta^2(s,z)).$$

It is obvious that $|\alpha'_{n,\varepsilon}(s,z)|$ is smaller than

$$K(\|X''(\varepsilon)_{s-} - X''(\varepsilon)_{(i-1)\Delta_n}\|^2 + \|X''(\varepsilon)_{s-} - X''(\varepsilon)_{(i-1)\Delta_n}\|^{m+l-2})\gamma(z)^2.$$

Therefore, (8.17) is a simple consequence of (8.14) and $\Gamma'(\varepsilon) \leq K$.  □

8.3. *Proof of Theorem 3.1*(a).  Observe that $B_T'^1 > 0$ and $B_T'^2 > 0$ on the set $\Omega_T^{(j)} \cup \Omega_T^{(d)}$, whereas $B_T > 0$ on $\Omega_T^{(j)}$ and $B_T = 0$ on $\Omega_T^{(d)}$. So, (3.16) is a trivial consequence of (4.3), which in turn comes from the following lemma.

LEMMA 8.3.  *If $h$ is a continuous function on $\mathbb{R}^2$ such that $h(x) = o(\|x\|^2)$ as $x \to 0$, we have $V(h, \Delta_n)_t \xrightarrow{\mathbb{P}} V(h)_t$ for each $t > 0$. We even have the convergence in probability (for the Skorokhod topology) of the processes $V(h, \Delta_n)$ toward $V(h)$.*

PROOF. Since $V(h)$ has no fixed times of discontinuity, the last claim implies the first one. When $h$ vanishes around the origin the result is proved exactly as in step 2 of Theorem 2.2 of [8]: the dimension of $X$ plays no role here.

Next, we turn to the general case. If $\varepsilon > 0$, we have $V(h(1-\psi_\varepsilon), \Delta_n) \xrightarrow{\mathbb{P}} V(h(1-\psi_\varepsilon))$ (for the Skorokhod topology) from what precedes, whereas $V(h(1-\psi_\varepsilon))$ obviously converges locally uniformly in time (for each $\omega$) to $V(h)_t$ as $\varepsilon \to 0$ by our assumptions on $h$, and hence, it is enough to prove that

$$(8.19) \quad \lim_{\varepsilon \to 0} \limsup_n \mathbb{P}\left( \sum_{s \le T} |V(h\psi_\varepsilon, \Delta_n)_t| > \eta \right) = 0 \qquad \forall \eta > 0, \ \forall T > 0.$$

We have $|h(x)| \le \theta(x_1) + \theta(x_2)$, where $\theta$ is a continuous function on $\mathbb{R}$ with $\theta(y) = o(y^2)$ as $y \to 0$. It is enough to prove (8.19) with $V(h\psi_\varepsilon, \Delta_n)_t$ substituted with $V^{(j)}(\theta\psi_\varepsilon, \Delta_n)_t = \sum_{i=1}^{[t/\Delta_n]} \theta(\Delta_i^n X^j) \psi(|\Delta_i^n X^j|/\varepsilon)$, for $j = 1, 2$. That is, we only need to prove (8.19) in the one-dimensional case, and this is a consequence of (3.4) in [8]. $\square$

8.4. *Proof of Theorem 4.1*(a). We start with a general result, of independent interest.

THEOREM 8.4. *Let $\phi$ be a $C^2$ function on $\mathbb{R}^2$ satisfying $\phi(0) = \phi_i'(0) = 0$ and $\phi_{ij}''(x) = o(\|x\|)$ as $x \to 0$, where $\phi_i'$ and $\phi_{ij}''$ are the first and second order partial derivatives. The two-dimensional processes*

$$(8.20) \quad \frac{1}{\sqrt{\Delta_n}}(V(\phi, \Delta_n)_t - V(\phi)_{\Delta_n[t/\Delta_n]}, V(\phi, k\Delta_n)_t - V(\phi)_{\Delta_n[t/k\Delta_n]})$$

*converge stably in law, on the product $\mathbb{D}(\mathbb{R}_+, \mathbb{R}) \times \mathbb{D}(\mathbb{R}_+, \mathbb{R})$ of the Skorokhod spaces, to the process with components*

$$(8.21) \quad (Z^1(\phi_1') + Z^2(\phi_2'), Z^1(\phi_1') + Z^2(\phi_2') + Z'^1(\phi_1') + Z'^2(\phi_2')).$$

We have the (stable) convergence in law of the processes in (8.20), as elements of the product functional space $\mathbb{D}(\mathbb{R}_+, \mathbb{R})^2$, but usually not as elements of the space $\mathbb{D}(\mathbb{R}_+, \mathbb{R}^2)$ with the (two-dimensional) Skorokhod topology, because a jump of $X$ entails a jump for both components above, but "with a probability close to $j/k$" the times at which these two components jump differ by an amount $j\Delta_n$, for $j = 1, \ldots, k-1$.

PROOF. The proof is essentially the same as for Theorem 2.12(i) of [8]. Fix $\varepsilon \in (0, \gamma_0]$, and let $S_q' = S_q'(\varepsilon)$ be the successive jump times of the Poisson process $\mu([0, t] \times A_\varepsilon^c)$, so that $X_t = X(\varepsilon)_t + \sum_{q: S_q' \le t} \Delta X_{S_q'}$ [notation (8.9)].

Next, we introduce some sets in which for this proof we could take $k_n = k$, but which are also needed later with $k_n$ as in (4.6). Namely, $\Omega_n(t, \varepsilon)$ denotes the set of all $\omega$ such that each interval $[0, t] \cap (i\Delta_n, (i + k_n)\Delta_n]$ contains at most one $S'_q$, and the intervals $(0, k_n\Delta_n)$ and $[t - (k_n + 1)\Delta_n, t]$ contains no $S'_q$, and finally $\|\Delta_i^n X(\varepsilon)\| \leq 2\varepsilon$ for all $i \leq t/\Delta_n$. Then,

$$(8.22) \qquad \Omega_n(t, \varepsilon) \to \Omega \qquad \text{as } n \to \infty \ \forall t, \varepsilon > 0.$$

We define the following variables on each set $\{(ik + j)\Delta_n < S'_q \leq (ik + j + 1)\Delta_n\}$ (with $0 \leq j < k$) as

$$\begin{cases} \bullet \ R_-(n, q) = X(\varepsilon)_{(ik+j)\Delta_n} - X(\varepsilon)_{ik\Delta_n}, \\ \bullet \ R_+(n, q) = X(\varepsilon)_{(i+1)k\Delta_n} - X(\varepsilon)_{(ik+j+1)\Delta_n}, \\ \bullet \ R_q^n = \Delta_{ik+j+1}^n X(\varepsilon), \qquad R_q'^n = R_-(n, q) + R_+(n, q), \qquad R_q''^n = R_q^n + R_q'^n. \end{cases}$$

Exactly as in [1], we have ($\xrightarrow{\mathcal{L}-(s)}$ denoting the stable convergence in law)

$$(8.23) \qquad (R_q^n/\sqrt{\Delta_n}, R_q'^n/\sqrt{\Delta_n})_{q \geq 1} \xrightarrow{\mathcal{L}-(s)} (R_q, R_q')_{q \geq 1}.$$

For any process $Y$, set $W(Y; \phi, \Delta_n)_t = V(Y; \phi, \Delta_n)_t - V(Y; \phi)_{\Delta_n[t/\Delta_n]}$. On the set $\Omega_n(T, \varepsilon)$, we have, for all $t \leq T$ and for $k' = 1$ and $k' = k$,

$$(8.24) \qquad W(X; \phi, k'\Delta_n)_t = W(X(\varepsilon); \phi, k'\Delta_n)_t + Y^{(\varepsilon)}(k'\Delta_n)_t,$$

where

$$Y^{(\varepsilon)}(\Delta_n)_t = \sum_{q \,:\, S'_q \leq \Delta_n[t/\Delta_n]} (\phi(\Delta X_{S'_q} + R_q^n) - \phi(\Delta X_{S'_q}) - \phi(R_q^n)),$$

$$Y^{(\varepsilon)}(k\Delta_n)_t = \sum_{q \,:\, S'_q \leq k\Delta_n[t/k\Delta_n]} (\phi(\Delta X_{S'_q} + R_q''^n) - \phi(\Delta X_{S'_q}) - \phi(R_q''^n)).$$

Since $\phi$ is $C^2$, we have

$$(8.25) \quad \begin{cases} Y^{(\varepsilon)}(\Delta_n)_t = \displaystyle\sum_{q \,:\, S'_q \leq \Delta_n[t/\Delta_n]} \left( \sum_{i=1}^2 \phi'_i(\Delta X_{S'_q} + \widetilde{R}_q^n) R^{n,i}_q - \phi(R_q^n) \right), \\ Y^{(\varepsilon)}(k\Delta_n)_t = \displaystyle\sum_{q \,:\, S'_q \leq k\Delta_n[t/k\Delta_n]} \left( \sum_{i=1}^2 \phi'_i(\Delta X_{S'_q} + \widetilde{R}_q''^n) R''^{n,i}_q - \phi(R_q''^n) \right), \end{cases}$$

where $\widetilde{R}_q^n$ and $\widetilde{R}_q''^n$ are between $0$ and $R_q^n$ and between $0$ and $\Delta R_q''^n$, respectively. Moreover, $\|\phi(x)\| = o(\|x\|^3)$; hence, by (8.23),

$$\frac{1}{\sqrt{\Delta_n}}(Y^{(\varepsilon)}(\Delta_n), Y^{(\varepsilon)}(k\Delta_n)) \xrightarrow{\mathcal{L}-(s)} \left( \sum_{i=1}^2 Z^{(\varepsilon)i}(\phi'_i), \sum_{i=1}^2 (Z^{(\varepsilon)i}(\phi_i) + Z'^{(\varepsilon)i}(\phi'_i)) \right)$$

[for the product topology of $\mathbb{D}(\mathbb{R}_+, \mathbb{R}) \times \mathbb{D}(\mathbb{R}_+, \mathbb{R})$], where $Z^{(\varepsilon)i}(\phi)$ and $Z'^{(\varepsilon)i}(\phi)$ are defined by (8.3), but with the sum taken only over the $S'_q(\varepsilon)$. As $\varepsilon \to 0$, the right-hand side above goes locally uniformly in time to the right-hand side of (8.21). Hence, in view of (8.24) and (8.22), it remains to prove that

$$(8.26) \qquad \lim_{\varepsilon \to 0} \limsup_n \mathbb{P}\left(\sup_{t \leq T} \frac{1}{\sqrt{\Delta_n}} |W(X(\varepsilon); \phi, k'\Delta_n)_t| > \eta\right) = 0$$

for all $\eta > 0$ and $T > 0$ and for $k' = 1$ and $k' = k$.

Now, with the notation (8.7), we set $\phi_\varepsilon(x) = \phi(x) \prod_{i=1}^2 \psi(|x_i|/2k\varepsilon)$. Then $\phi_\varepsilon$ is a $C^2$ function which coincides with $\phi$ when $\|x\| \leq 2k\varepsilon$ and vanishes for $\|x\| > \sqrt{2} \times 4k\varepsilon$. Hence for each $T$, on a set of probability going to 1 as $n \to \infty$ we have $W(X(\varepsilon); \phi, k'\Delta_n)_t = W(X(\varepsilon); \phi_\varepsilon, k'\Delta_n)_t$ for all $t \leq T$. Therefore it is enough to prove (8.26) with $\phi_\varepsilon$ instead of $\phi$. But this is exactly the last step in the proof of Theorem 2-11(i) of [8] (in which the dimension of $X$ plays no role; this is where the hypothesis $\phi''_{ij}(x) = o(\|x\|)$ is used). Hence, we are done. $\square$

LEMMA 8.5. *For any real-valued continuous function $\phi$ on $\mathbb{R}^2$ such that $\phi(x) = O(\|x\|^2)$ as $x \to 0$, we have*

$$(8.27) \qquad \frac{1}{\sqrt{\Delta_n}} \sum_{k\Delta_n[T/k\Delta_n] < s \leq T} |\phi(\Delta X_s)| \xrightarrow{\mathbb{P}} 0.$$

PROOF. Denote, by $U_n$, the left-hand side of (8.27). Assumption (SH) yields $|\phi(\delta(s,z))| \leq K\gamma(z)^2$ for some constant $K$ (recall that here $\delta$ is bounded); hence,

$$\mathbb{E}(U_n) \leq \frac{K}{\sqrt{\Delta_n}} \mathbb{E}\left(\int_{k\Delta_n[T/k\Delta_n]}^T ds \int \gamma(z)^2 \lambda(dz)\right) \leq K' k \sqrt{\Delta_n}$$

for another constant $K'$, and the result follows. $\square$

PROOF OF THEOREM 4.1(a). We apply Theorem 8.4 with $\phi = f$. Since $X$ has no fixed time of discontinuity, and by the previous lemma, we deduce

$$G_T^n := \frac{1}{\sqrt{\Delta_n}}(V(f, k\Delta_n)_T - V(f, \Delta_n)_T) \xrightarrow{\mathcal{L}-(s)} \widetilde{G}_T.$$

We also have $\Phi_n^{(j)} - 1 = \sqrt{\Delta_n} G_T^n / V(f, \Delta_n)_T$; hence, $(\Phi_n^{(j)} - 1)\sqrt{\Delta_n}$ converges stably in law, in restriction to the set $\Omega_T^{(j)} = \{B_T > 0\}$, to $\widetilde{G}_T / B_T$ by Lemma 8.3. The end of the claim follows from Lemma 8.1 and from (8.5). $\square$

8.5. *Proof of Theorems 3.1(b) and 4.1(b).* Equation (3.17) follows from Lemma 8.3. For (3.18) and Theorem 4.1(b), we will use the following theorem.

THEOREM 8.6. *Let $h$ be a $d$-dimensional $C^2$ function on $\mathbb{R}^2$, its first component being $f$, and all the others being either vanishing on a neighborhood of $0$ or equal to $x \mapsto x_1^m x_2^l$ for some $m, l \geq 2$ with $m + l \geq 5$. Assume also that (with $h_i'$ and $h_{ij}''$ being the $\mathbb{R}^d$-valued partial derivatives)*

$$(8.28) \qquad \begin{aligned} x &= (x_1, 0) \quad or \\ x &= (0, x_2) \quad \Rightarrow \quad h(x) = h_1'(x) = h_2'(x) = h_{12}''(x) = 0. \end{aligned}$$

*Then, the $2d$-dimensional processes*

$$(8.29) \qquad \left( \frac{1}{\Delta_n} V(h, \Delta_n)_t, \frac{1}{\Delta_n} V(h, k\Delta_n)_t \right)_{t \in [0, T]}$$

*converge stably in law, in restriction to the union $\Omega_T^{(d)} \cup \Omega_T^{(c)}$ and on the product $\mathbb{D}([0, T], \mathbb{R}^d) \times \mathbb{D}([0, T], \mathbb{R}^d)$, to the process $(\widetilde{D}(h)_t + C(h)_t, \widetilde{D}''(h)_t + kC(h)_t)_{t \in [0, T]}$, where*

$$(8.30) \qquad \begin{cases} \widetilde{D}(h)_t = \frac{1}{2} \sum\limits_{q:\, S_q \leq t} (h_{11}''(\Delta X_{S_q})(R_q^1)^2 + h_{22}''(\Delta X_{S_q})(R_q^2)^2), \\ \widetilde{D}''(h)_t = \frac{1}{2} \sum\limits_{q:\, S_q \leq t} (h_{11}''(\Delta X_{S_q})(R_q'^1)^2 + h_{22}''(\Delta X_{S_q})(R_q''^2)^2) \end{cases}$$

*and where $C(h)$ is the process whose first component is $C$, as given by (3.12), and all others are $0$.*

Any component of the form $x_1^m x_2^l$ satisfies (8.28), so this condition is a condition on the components which vanish on a neighborhood of $0$.

PROOF OF THEOREM 8.6. (1) As said before, we assume Assumption (SH). But another localization allows to do more: let $\tau_q = \inf(t : \|\delta_t'\| \geq q)$ [the process $\delta'$ is defined in (2.3)]. By (d) of (H) we have $\lim_q \tau_q \geq \tau$, and of course $\tau > T$ on the set $\Omega_T^{(d)} \cup \Omega_T^{(c)}$, so it is enough to prove the result in restriction to each set $(\Omega_T^{(d)} \cup \Omega_T^{(c)}) \cap \{\tau_q > T\}$. Now, let $X^{(q)}$ be the process defined by (2.1), with the coefficients

$$b_t^{(q)} = (b_t - \delta_t') 1_{\{t \leq \tau_q\}}, \qquad \sigma_t^{(q)} = \sigma_t, \qquad \delta^{(q)}(t, z) = \delta(t, z) 1_{\widetilde{\Gamma}^c}(t, z).$$

These coefficients satisfy Assumption (SH). Moreover, we obviously have $X_t = X_t^{(q)}$ for all $t \leq T$ on the set $(\Omega_T^{(j)})^c \cap \{\tau_q > T\}$; hence, the process (8.29) and also $\widetilde{D}(h)_t$, $\widetilde{D}''(h)_t$ and $C(h)_t$ for $t \leq T$ are the same on that set, whether computed on the basis of $X$ or on the basis of $X^{(q)}$. This means that, for proving our result, we can substitute $X$ with $X^{(q)}$, which satisfies Assumption (SH) and whose two components have no common jumps by construction.

In other words, we can and will assume in the rest of the proof that $X$ satisfies Assumption (SH) and that $X^1$ and $X^2$ have no common jumps. We will then prove that, in fact, the stable convergence in law holds everywhere [not only on $(\Omega_T^{(j)})^c$], and for the time interval $\mathbb{R}_+$. We use the same notation as in the proof of Theorem 8.4.

(2) Pick $\varepsilon \in (0, \gamma_0]$, and take $S'_q = S'_q(\varepsilon)$. Exactly as for (8.24), on the set $\Omega_n(T, \varepsilon)$ and for $t \leq T$ and $k' = 1$ or $k' = k$, we have

$$(8.31) \qquad V(X; h, k'\Delta_n)_t = V(X(\varepsilon); h, k'\Delta_n)_t + \overline{Y}^{(\varepsilon)}(h, k'\Delta_n)_t,$$

where

$$\overline{Y}^{(\varepsilon)}(h, \Delta_n)_t = \sum_{q : S'_q \leq \Delta_n[t/\Delta_n]} (h(\Delta X_{S'_q} + R_q^n) - h(R_q^n))$$

$$\overline{Y}^{(\varepsilon)}(h, k\Delta_n)_t = \sum_{q : S'_q \leq k\Delta_n[t/k\Delta_n]} (h(\Delta X_{S'_q} + R_q''^n) - h(R_q''^n)).$$

By a Taylor expansion and the properties $h(\Delta X_{S'_q}) = h'_i(\Delta X_{S'_q}) = 0$, those are, respectively, equal to

$$\sum_{q : S'_q \leq \Delta_n[t/\Delta_n]} \left( \tfrac{1}{2} \sum_{i,j=1}^{2} h''_{ij}(\Delta X_{S'_q} + \widetilde{R}_q^n) R_q^{n,i} R_q^{n,j} - h(R_q^n) \right),$$

$$\sum_{q : S'_q \leq k\Delta_n[t/k\Delta_n]} \left( \tfrac{1}{2} \sum_{i,j=1}^{2} h''_{ij}(\Delta X_{S'_q} + \widetilde{R}_q''^n) R_q''^{n,i} R_q''^{n,j} - h(R_q''^n) \right),$$

where $\widetilde{R}_q^n$, respectively, $\widetilde{R}_q''^n$, is between 0 and $R_q^n$, respectively, $R_q''^n$. Since we also have $h''_{12}(\Delta X_{S'_q}) = 0$, (8.23) yields

$$(\Delta_n^{-1}\overline{Y}^{(\varepsilon)}(h, \Delta_n), \Delta_n^{-1}\overline{Y}^{(\varepsilon)}(h, k\Delta_n)) \xrightarrow{\mathcal{L}-(s)} (\widetilde{D}^{(\varepsilon)}(h), \widetilde{D}''^{(\varepsilon)}(h))$$

[for the product topology of $\mathbb{D}(\mathbb{R}_+, \mathbb{R}^d) \times \mathbb{D}(\mathbb{R}_+, \mathbb{R}^d)$], where $\widetilde{D}^{(\varepsilon)}(h)$ and $\widetilde{D}''^{(\varepsilon)}(h)$ are defined by (8.30), but with the sum taken over the $S'_q(\varepsilon)$ only. As $\varepsilon \to 0$, we have $\widetilde{D}^{(\varepsilon)}(h)_t \to \widetilde{D}(h)_t$ and $\widetilde{D}''^{(\varepsilon)}(h)_t \to \widetilde{D}''(h)_t$ locally uniformly in $t$. Hence, in view of (8.31) and (8.22) it remains to prove that, for all $\eta > 0$ and $k' = 1$ and $k' = k$,

$$(8.32) \quad \lim_{\varepsilon \to 0} \limsup_n \mathbb{P}\left( \sup_{t \leq T} \left\| \frac{1}{\Delta_n} V(X(\varepsilon); h, k'\Delta_n)_t - k'C(h)_t \right\| > \eta \right) = 0.$$

(3) Obviously, it suffices to show (8.32) for each component or, equivalently, we can assume that $h$ is one-dimensional. If $h(x) = 0$ when $\|x\| \leq \rho$, then, since $X(\varepsilon)$ has jumps smaller than $\varepsilon$, we see that if $\varepsilon < \rho/2$, $V(X(\varepsilon); h, k'\Delta_n)_t$ vanishes for all $t \in [0, T]$ on the set $\Omega_n(T, \varepsilon)$ [see (8.22)].

Therefore, in this case, (8.32) is obvious, and it remains to study the case where $h(x) = x_1^m x_2^l$ for $m, l \geq 2$.

(4) Recall that $X(\varepsilon) = X'(\varepsilon) + X''(\varepsilon)$ [see (8.10)]. The process $X'(\varepsilon)$ is a continuous Itô semimartingale with bounded coefficients and càdlàg volatility. Since $h$ is homogeneous of degree $r = m + l$ it is known (see, e.g., Theorem 2.4(i) of [8] this theorem is for a one-dimensional process, but the multidimensional extension is straightforward; also, see [11]) that for each fixed $\varepsilon > 0$, the processes $\Delta_n^{1-r/2} V(X'(\varepsilon); h, k'\Delta_n)$ converge in probability, locally uniformly in time, to a limit, that is, $k'C$ when $r = 4$ (i.e., when $h$ is the function $f$) and that need not be specified when $r \geq 5$. Then, $\Delta_n^{-1} V(X'(\varepsilon); h, k'\Delta_n)$ converges to $k'C(h)$, in probability, locally uniformly in time, and, in order to obtain (8.32), it is clearly enough to show that

$$\lim_{\varepsilon \to 0} \limsup_n \mathbb{P}\left( \sup_{t \leq T} \frac{1}{\Delta_n} |V(X(\varepsilon); h, k'\Delta_n)_t - V(X'(\varepsilon), h, k'\Delta_n)_t| > \eta \right)$$
$$= 0. \tag{8.33}$$

(5) We prove (8.33) for $k' = 1$, the proof for $k' = k$ being similar. For all $u > 0$, $v, w \geq 0$ and $p, q \geq 1$ we have $v^p w^q \leq u v^{p+q} + w^{p+q} u^{-p/q}$, hence since $h(x) = (x_1)^m (x_2)^l$ we see that for all $u > 0$ there is a constant $A_u$ (depending also on $m, l$) such that, for all $x, y \in \mathbb{R}^2$,

$$|h(x + y) - h(x)| \leq u|h(x)| + A_u(|x_1|^m |y_2|^l + |y_1|^m |x_2|^l + |y_1|^m |y_1|^l).$$

It follows that

$$|V(X(\varepsilon); h, \Delta_n) - V(X'(\varepsilon); h, \Delta_n)| \leq u V(X'(\varepsilon); |h|, \Delta_n) + A_u U^n(\varepsilon), \tag{8.34}$$

where $U^n(\varepsilon)_t = \sum_{i=1}^{[t/\Delta_n]} \zeta_i^n(\varepsilon)$ and

$$\zeta_i^n(\varepsilon) = |\Delta_i^n X''^1(\varepsilon)|^m |\Delta_i^n X''^2(\varepsilon)|^l + |\Delta_i^n X''^1(\varepsilon)|^m |\Delta_i^n X'^2(\varepsilon)|^l$$
$$+ |\Delta_i^n X'^1(\varepsilon)|^m |\Delta_i^n X''^2(\varepsilon)|^l.$$

First, (8.16) and (8.17) yield $\mathbb{E}(\zeta_i^n(\varepsilon)) \leq K\sqrt{\Gamma'(\varepsilon)}\Delta_n^2(1 + \Delta_n^{1/4}/\varepsilon^{m \vee l})$. Therefore, since $\Gamma'(\varepsilon) \to 0$ as $\varepsilon \to 0$, we obtain

$$\lim_{\varepsilon \to 0} \limsup_n \mathbb{E}\left( \frac{1}{\Delta_n} U^n(\varepsilon)_T \right) = 0.$$

Second, as seen above, $Z_n := \Delta_n^{-1} V(X'(\varepsilon); |h|, \Delta_n)_T - C(h)_T \xrightarrow{\mathbb{P}} 0$; hence,

$$\lim_{u \to 0} \limsup_n \mathbb{P}\left( \frac{u}{\Delta_n} V(X'(\varepsilon); |h|, \Delta_n)_T > \frac{\eta}{2} \right)$$
$$\leq \lim_{u \to 0} \limsup_n \left( \mathbb{P}\left( |Z_n| > \frac{\eta}{4u} \right) + \mathbb{P}\left( C(h)_T > \frac{\eta}{4u} \right) \right) = 0.$$

These two results, put together with (8.34), allow us to deduce (8.33). □

PROOF OF THEOREM 3.1(b). First, the variable $\widetilde{\Phi}$ defined by (3.18) on $\Omega_T^{(d)}$ and, say 1 elsewhere, takes its values in $(0, \infty)$, and conditionally on $\mathcal{F}$, and, in restriction to $\Omega_T^{(d)}$, it has a density [recall (3.11) and (3.14)], so $\widetilde{\Phi} \neq 1$ a.s. on $\Omega_T^{(d)}$. Second, the convergence $\Phi_n^{(j)} \overset{\mathcal{L}-(s)}{\longrightarrow} \widetilde{\Phi}$, in restriction to $\Omega_T^{(d)}$ is obvious from Theorem 8.6. □

PROOF OF THEOREM 4.1(b). By Lemma 8.3 applied to $h = g_1$ and $h = g_2$, and by Theorem 8.6, it is obvious that $\Phi_n^{(d)}/\Delta_n \overset{\mathcal{L}-(s)}{\longrightarrow} \widetilde{\Phi}' = (\widetilde{D}_T + C_T)/\sqrt{B_T'^1 B_T'^2}$, in restriction to $\Omega_T^{(d)}$. Finally, (4.2) follows from (8.1). □

8.6. *Proof of (4.10).* The first part of (4.10) is none other than a multi-dimensional version of (26) in [1], applied with $q = 4$ and $r = 1$, and we leave the (simple) computations to the reader. As to the second part, it is almost trivial. Indeed, setting $f_\rho = f \psi_\rho$ [notation (8.7)], for any $\rho > 0$ we have

$$\Delta_n \widehat{A}'(\Delta_n)_T \leq V(f_\rho, \Delta_n)_T$$

as soon as $\alpha \Delta_n^\varpi < \rho$. Now, Lemma 8.3 yields that $V(f_\rho, \Delta_n)_T$ converges in probability to $V(f_\rho)_T$, which in turn goes to 0 as $\rho \to 0$; hence, we have the result.

8.7. *Proof of (4.11).* Both claims in (4.11) amount to prove the following property: introduce the functions $g(x) = x_1^u x_2^v$ where $u + v \geq 2$, and $\overline{g}(x) = x_m x_l$ where $m, l$ are two indices taking the values 1 or 2. We complete this notation with $g_n(x) = g(x) 1_{\{\|x\| > \alpha \Delta_n^\varpi\}}$ and $g_{n,\rho} = g_n \psi_\rho$ and $g'_{n,\rho} = g_n - g_{n,\rho}$ for $\rho > 0$, and $\overline{g}_n(x) = \overline{g}(x) 1_{\{\|x\| \leq \alpha \Delta_n^\varpi\}}$. Then, we need to prove that

$$\begin{align}
(8.35) \quad \widehat{H}_t^n &:= \frac{1}{k_n \Delta_n} \sum_{i=1+k_n}^{[t/\Delta_n]-k_n-1} g_n(\Delta_i^n X) \sum_{j \in I_n(i)} \overline{g}_n(\Delta_j^n X) \\
&\overset{\mathbb{P}}{\longrightarrow} H_t = \sum_{s \leq t} g(\Delta X_s)(c_{s-}^{ml} + c_s^{ml}).
\end{align}$$

The proof is basically the same as for (27) of [1], and goes through several steps.

*Step* 1. This step is devoted to showing some estimates. Recall (8.11). First, $|g_{n,\rho}(x)| \leq \|x\|^{u+v}$ and $|g'_{n,\rho}(x)| \leq \|x\|^{u+v}$ with $u+v \geq 2$, and $|\overline{g}_n(x)| \leq \|x\|^2$; hence, (8.12) and (8.14) yield

$$(8.36) \quad \mathbb{E}((|g_{n,\rho}| + |g'_{n,\rho}| + |\overline{g}_n|)(\Delta_i^n X) \mid \mathcal{F}_{(i-1)\Delta_n}) \leq K\Delta_n.$$

Second, when $\rho \leq 1/2$, we have

$$|g_{n,\rho}(x+y)| \leq K(\|x\|^4 \Delta_n^{-2\varpi} + (\|y\|^2 \wedge \rho^2)).$$

Using this with $x = \Delta_i^n X'$ and $y = \Delta_i^n X''$, plus (8.12) and (8.15) with $\eta = \rho$ and $\theta = \sqrt{\rho}$, we obtain

$$(8.37) \quad \mathbb{E}(|g_{n,\rho}(\Delta_i^n X)| \mid \mathcal{F}_{(i-1)\Delta_n}) \leq \Delta_n a_n(\rho), \qquad \lim_{\rho \to 0} \lim_{n \to \infty} a_n(\rho) = 0,$$

provided we take $a_n(\rho) = K(\Delta_n^{1-2\varpi} + \rho + \frac{\Delta_n}{\rho} + \Gamma''(\rho))$.

Next, we set $\delta_j^n = \sigma_{(j-1)\Delta_n} \Delta_j^n W$. It is easily checked that for all $w > 0$ there is a constant $A_w$ such that for $x, y, z \in \mathbb{R}^2$,

$$|\overline{g}_n(x+y+z) - \overline{g}(x)| \leq w\|x\|^2 + A_w \left( \frac{\|x\|^4 + \|y\|^4}{\Delta_n^{2\varpi}} + \|y\|^2 + (\|z\|^2 \wedge \Delta_n^{2\varpi}) \right).$$

If we apply this with $x = \delta_j^n$ and $y = \Delta_j^n X' - \delta_j^n$ and $z = \Delta_j^n X''$, plus (8.13) and (8.15) with $\eta = \theta^2 = \Delta_n^\varpi$, we obtain, with $a_n' = \Delta_n^{\varpi \wedge (1-2\varpi)} + \Gamma''(\Delta_n^{\varpi/2})$,

$$\mathbb{E}(|\overline{g}_n(\Delta_j^n X) - \overline{g}(\delta_j^n)| \mid \mathcal{F}_{(j-1)\Delta_n})$$
$$(8.38) \qquad \leq K\Delta_n w + K_w \Delta_n a_n' + K_w Y_j^n,$$
$$Y_j^n = \mathbb{E}\left( \int_{(j-1)\Delta_n}^{j\Delta_n} \|\sigma_s - \sigma_{(j-1)\Delta_n}\|^2 \, ds \mid \mathcal{F}_{(j-1)\Delta_n} \right).$$

*Step* 2. For each $\rho > 0$, we set

$$\widehat{H}(\rho)_t^n = \frac{1}{k_n \Delta_n} \sum_{i=1+k_n}^{[t/\Delta_n]-k_n-1} g_{n,\rho}(\Delta_i^n X) \sum_{j \in I_n(i)} \overline{g}_n(\Delta_j^n X),$$

$$\widehat{H}(\rho)_t'^n = \frac{1}{k_n \Delta_n} \sum_{i=1+k_n}^{[t/\Delta_n]-k_n-1} g_{n,\rho}'(\Delta_i^n X) \sum_{j \in I_n(i)} (\overline{g}_n(\Delta_j^n X) - \overline{g}(\delta_j^n)),$$

$$\overline{H}(\rho)_t^n = \frac{1}{k_n \Delta_n} \sum_{i=1+k_n}^{[t/\Delta_n]-k_n-1} g_{n,\rho}'(\Delta_i^n X) \sum_{j \in I_n(i)} \overline{g}(\delta_j^n),$$

$$H(\rho)_t = \sum_{q : S_q \leq t} (g\psi_\rho)(\Delta X_{S_q})(c_{S_q-}^{ml} + c_{S_q}^{ml}), \qquad \overline{H}(\rho) = H - H(\rho).$$

We have $\widehat{H}^n = \widehat{H}(\rho)^n + \widehat{H}(\rho)'^n + \overline{H}(\rho)^n$; hence, for (8.35), it is enough to prove the following four properties:

$$(8.39) \qquad\qquad \rho \to 0 \quad \Rightarrow \quad H(\rho)_t \xrightarrow{\mathbb{P}} 0,$$

$$(8.40) \qquad \lim_{\rho \to 0} \limsup_n \mathbb{E}(|\widehat{H}(\rho)_t^n|) = 0,$$

$$(8.41) \qquad \lim_{\rho \to 0} \limsup_n \mathbb{E}(|\widehat{H}(\rho)_t'^n|) = 0,$$

$$(8.42) \qquad \rho \in (0,1), \qquad n \to \infty \quad \Rightarrow \quad \overline{H}(\rho)_t^n \xrightarrow{\mathbb{P}} \overline{H}(\rho)_t.$$

Note that the property (8.39) readily follows from Lebesgue theorem.

*Step* 3. Here, we prove (8.40) and (8.41). First, by successive conditioning, we deduce from (8.36) and (8.37) that

$$j \neq i \quad \Rightarrow \quad \mathbb{E}(|g_{n,\rho}(\Delta_i^n X) \overline{g}_n(\Delta_j^n X)|) \leq K \Delta_n^2 a_n(\rho).$$

Since $\lim_{\rho \to 0} \lim_n a_n(\rho) = 0$, we readily deduce (8.40). Second, again by successive conditioning, we deduce from (8.36) and (8.38) that, for all $w > 0$,

$$j \neq i \quad \Rightarrow \quad \mathbb{E}(|g_{n,\rho}'(\Delta_i^n X)(\overline{g}_n(\Delta_j^n X) - \overline{g}(\delta_j^n))|)$$
$$\leq K \Delta_n^2 w + K_w \Delta_n^2 a_n' + K_w \Delta_n \mathbb{E}(Y_j^n).$$

This readily yields

$$\mathbb{E}(|\widehat{H}(\rho)_t'^n|) \leq K t a_n(\rho) + K t w + K_w t a_n' + K_w \mathbb{E}\left( \int_0^t \|\sigma_s - \sigma_{\Delta_n[s/\Delta_n]}\|^2 \, ds \right).$$

Observe that $\|\sigma_s - \sigma_{\Delta_n[s/\Delta_n]}\|$ is bounded [uniformly in $(\omega, s, n)$], and goes to 0 for $\mathbb{P}(d\omega) \otimes ds$ almost all $(\omega, s)$ as $n \to \infty$ because $\sigma_s$ is right continuous with left limits. Then, by Lebesgue convergence theorem and $a_n' \to 0$, we get

$$\limsup_n \mathbb{E}(|\widehat{H}(\rho)_t'^n|) \leq K t \lim_n a_n(\rho) + K t w$$

and since $w > 0$ is arbitrary and $\lim_{\rho \to 0} \lim_n a_n(\rho) = 0$, we deduce (8.41).

*Step* 4. It remains to prove (8.42). Fix $\rho \in (0,1)$ and $t > 0$, and recall the jump times $S_q' = S_q'(\rho/2)$ and the set $\Omega_n(T, \rho/2)$ of the proof of Theorem 8.4. On $\Omega_n(t, \rho/2)$ there is no $S_q'$ in $(0, k_n \Delta_n]$, nor in $(t - (k_n + 1)\Delta_n, t]$ and there is at most one $S_q'$ in an interval $((i-1)\Delta_n, i\Delta_n]$ with $i\Delta_n \leq t$, and if $((i-1)\Delta_n, i\Delta_n]$ contains no $S_q'$ we have $\psi_\rho(\Delta_i^n X) = 1$. Hence, on $\Omega_n(t, \rho/2)$,

$$\overline{H}(\rho)_t^n = \sum_{q : k_n \Delta_n < S_q' \leq t - (k_n+1)\Delta_n} g_{n,\rho}'(\Delta_{i(n,q)}^n X) \frac{1}{k_n \Delta_n} \sum_{j \in I_n(i(n,q))} \overline{g}(\delta_j^n),$$

where $i(n,q) = \inf(i : i\Delta_n \geq S_q')$. Observe also that

$$\overline{H}(\rho)_t = \sum_{q : S_q' \leq t} (g(1 - \psi_\rho))(\Delta X_{S_q'})(c_{S_q'}^{ml} + c_{S_q'}^{ml}).$$

The sum over $q$ with $S_q' \leq t$ is finite, and obviously $g_{n,\rho}'(\Delta_{i(n,q)}^n X) \to (g(1 - \psi_\rho))(\Delta X_{S_q'})$ pointwise. Hence, for (8.42), and since $\Omega_n(t, \rho/2) \to \Omega$ as $n \to$

$\infty$, we need only to prove that

$$
\text{(8.43)}
\begin{aligned}
\frac{1}{k_n \Delta_n} \sum_{j \in I_{n,-}(i(n,q))} \overline{g}(\delta_j^n) &\xrightarrow{\mathbb{P}} c_{S_q'-}, \\
\frac{1}{k_n \Delta_n} \sum_{j \in I_{n,+}(i(n,q))} \overline{g}(\delta_j^n) &\xrightarrow{\mathbb{P}} c_{S_q'}.
\end{aligned}
$$

This is proved in (71) of [1] when $X$ is one-dimensional, and the two-dimensional extension is straightforward.

At this point, the proof of (4.11) is finished. However, we will now derive a consequence of (8.43), to be used later.

LEMMA 8.7. *In the previous setting, and in particular with $\rho > 0$ fixed, for any $q \geq 1$, we have*

$$
\text{(8.44)} \qquad \widehat{c}(n,-)_{i(n,q)} \xrightarrow{\mathbb{P}} c_{S_q'-}, \qquad \widehat{c}(n,+)_{i(n,q)} \xrightarrow{\mathbb{P}} c_{S_q'}.
$$

PROOF. In view of (8.43), this is a simple consequence of

$$
\text{(8.45)} \qquad U(n,q,\pm) := \frac{1}{k_n \Delta_n} \sum_{j \in I_{n,\pm}(i(n,q))} |\overline{g}_n(\Delta_j^n X) - \overline{g}(\delta_j^n)| \xrightarrow{\mathbb{P}} 0.
$$

We have a problem here: we cannot apply (8.38) without care, because the integer $i(n,q)$ is random, as it is a function of $S_q'$. We set $\zeta(n,q,+) = \sup(\|\sigma_s - \sigma_{S_q'}\|^2 : s \in [S_q', S_q' + k_n \Delta_n])$ and $\zeta(n,q,-) = \sup(\|\sigma_s - \sigma_{S_q'-}\|^2 : s \in [S_q' - k_n \Delta_n, S_q'))$, which are bounded and converge to 0 pointwise as $n \to \infty$, and $Y_j^n \leq \Delta_n \zeta(n,q,\pm)$ when $j \in I_{n,\pm}(i(n,q))$.

For $U(n,q,+)$, (8.45) is easy. Indeed, if $2 < j \leq k_n + 1$, we deduce from (8.38) and from the property $\{i(n,q) = r\} \in \mathcal{F}_{(r+1)\Delta_n}$ that

$$
\text{(8.46)}
\begin{aligned}
&\mathbb{E}(|\overline{g}_n(\Delta_{i(n,q)+j}^n X) - \overline{g}(\delta_{i(n,q)+j}^n)|) \\
&\quad = \sum_{r \geq 1} \mathbb{E}(|\overline{g}_n(\Delta_{r+j}^n X) - \overline{g}(\delta_{r+j}^n)| 1_{\{i(n,q)=r\}}) \\
&\quad = \sum_{r \geq 1} \mathbb{E}(\mathbb{E}(|\overline{g}_n(\Delta_{r+j}^n X) - \overline{g}(\delta_{r+j}^n)| \mid \mathcal{F}_{(r+1)\Delta_n}) 1_{\{i(n,q)=r\}}) \\
&\quad \leq \Delta_n \sum_{r \geq 1} \mathbb{E}((Kw + K_w a_n' + K_w \zeta(n,q,+)) 1_{\{i(n,q)=r\}}) \\
&\quad = \Delta_n (Kw + K_w a_n' + K_w \mathbb{E}(\zeta(n,q,+))).
\end{aligned}
$$

This holds for all $w > 0$, and $\mathbb{E}(\zeta(n,q,+)) \to 0$, we have (8.45).

For $U(n,q,-)$ things are more difficult. We replace the Poisson measure $\mu$ by its restriction $\overline{\mu}$ to $\mathbb{R}_+ \times A_{\rho/2}$, which is again Poisson, and we define $\overline{X}$ by (2.1), with $\mu$ substituted with $\overline{\mu}$, and we associate $\overline{\zeta}(n, q, -)$ and $\overline{U}(n, q, -)$ as above. We consider the filtration $\overline{\mathcal{F}}_t = \mathcal{F}_t \vee \sigma(\mu([0, s] \times A) : s \geq 0, A \in \mathcal{E}, A \cap A_{\rho/2} = \varnothing)$. Then, due to the independence properties of the jumps of the Poisson measure $\mu$ and of $W$ and $\mu$, $\overline{X}$ is again an Itô semimartingale relative to the filtration $(\overline{\mathcal{F}}_t)$, and thus all estimates for $X$, relative to the filtration $(\mathcal{F}_t)$, also hold for $\overline{X}$, relative to the filtration $(\overline{\mathcal{F}}_t)$.

Now, the random variable $i(n, q)$ becomes $\overline{\mathcal{F}}_0$-measurable, so we can argue as in (8.46) for the process $\overline{X}$, using in particular (8.38) with $\overline{\mathcal{F}}_{(j-1)\Delta_n}$, even when $j$ is negative, between $-k_n$ and $-1$, provided we add the indicator function of the set $\{i(n, p-1) + k_n < i(n, p)\}$. Hence, we deduce as before that $\mathbb{E}(\overline{U}(n, q, -)) \to 0$. It remains to observe that on the set $\Omega_n(t, \rho/2)$ we have $U(n, q, -) = \overline{U}(n, q, -)$ as soon as $S_q' \leq t$, and (8.45) for $U(n, q, -)$ is deduced from (8.22) and from the above. $\square$

8.8. *Proof of Theorem 4.2*(b). For (b) Theorem 4.2, and exactly as in the proof of Theorem 8.6, we can assume Assumption (SH) *and also* that $X^1$ and $X^2$ do not jump together, and, under this additional assumption, we prove following two convergences on $\Omega$:

$$(8.47) \qquad \widehat{A}_t'^n \overset{\mathbb{P}}{\longrightarrow} C_t,$$

$$(8.48) \qquad \text{the sequence of variables } \left(\frac{1}{\Delta_n} \widehat{F}_t'^n\right)_{n \geq 1} \text{ is tight.}$$

We begin with (8.48), which is easy. Indeed, it suffices to prove that $\mathbb{E}(\widehat{F}_t'^n) \leq Kt\Delta_n$ for some constant $K$. In view of (4.9), this amounts to proving that

$$(8.49) \qquad \mathbb{E}((\Delta_i^n X^m)^4 (\Delta_i^n X^l)^2 (\Delta_j^n X^l)^2) \leq K\Delta_n^3,$$

if $(m, l)$ equals $(1, 2)$ or $(2, 1)$ and $i \neq j$. Since $X^1$ and $X^2$ have no common jumps, (8.12) and (8.17) applied with $\varepsilon = \gamma_0$ yield that

$$\mathbb{E}((\Delta_i^n X^m)^4 (\Delta_i^n X^l)^2 \mid \mathcal{F}_{(i-1)\Delta_n}) \leq K\Delta_n^2, \qquad \mathbb{E}((\Delta_j^n X^l)^2 \mid \mathcal{F}_{(j-1)\Delta_n}) \leq K\Delta_n$$

and by successive integration we get (8.49).

The proof of (8.47) necessitates several steps.

*Step* 1. Observe that $f$, of course, but also $f_\rho = f\psi_\rho$ for any $\rho > 0$, satisfies (8.28). Then, by Theorem 8.6, the pairs $(\frac{1}{\Delta_n} V(f, \Delta_n), \frac{1}{\Delta_n} V(f - f_\rho, \Delta_n))$ converge stably in law to $(\widetilde{D}(f) + C, \widetilde{D}(f - f_\rho))$, for the Skorokhod topology on $\mathbb{D}(\mathbb{R}_+, \mathbb{R}^2)$ (since we assume that $\Omega_T^{(j)} = \varnothing$ for all $T$). Taking the difference of the two components, we deduce that, for any $\rho > 0$,

$$(8.50) \qquad \frac{1}{\Delta_n} V(f_\rho, \Delta_n) \overset{\mathcal{L}-(s)}{\longrightarrow} \widetilde{D}(f_\rho) + C.$$

*Step* 2. Now, we prove that, for any $A > 1$ and $t > 0$,

$$(8.51) \qquad \frac{1}{\Delta_n} V(f_{A\sqrt{\Delta_n}}, \Delta_n)_t \xrightarrow{\mathbb{P}} \int_0^t \rho_{\sigma_s}(f_A) \, ds,$$

where for any $2 \times 2$ matrix $\sigma$ we have set $\rho_\sigma(f_A) = \mathbb{E}(f_A(\sigma U))$, where $U$ is an $\mathcal{N}(0, I_2)$ variable [hence, $\rho_{\sigma_s}(f_A) = \rho_{\sigma_s(\omega)}(f_A)$ is a process].

Recalling (8.11), and by Theorem 2.4(i) of [8] (see [11] for the multidimensional version),

$$\frac{1}{\Delta_n} V(X'; f_{A\sqrt{\Delta_n}}, \Delta_n)_t \xrightarrow{\mathbb{P}} \int_0^t \rho_{\sigma_s}(f_A) \, ds;$$

hence, it is enough to prove that

$$(8.52) \quad G_n := \frac{1}{\Delta_n} |V(X; f_{A\sqrt{\Delta_n}}, \Delta_n)_t - V(X'; f_{A\sqrt{\Delta_n}}, \Delta_n)_t| \xrightarrow{\mathbb{P}} 0.$$

It is obvious that $|f_\eta(x + y) - f_\eta(x)| \leq K \eta^3 (\|y\| \wedge \eta)$ for all $x, y \in \mathbb{R}^2$, $\eta > 0$. Then, if we use (8.15) for $\varepsilon = \gamma_0$, we see that

$$\mathbb{E}(G_n) \leq \frac{1}{\Delta_n} \sum_{i=1}^{[t/\Delta_n]} \mathbb{E}(|f_{A\sqrt{\Delta_n}}(\Delta_i^n X' + \Delta_i^n X'') - f_{A\sqrt{\Delta_n}}(\Delta_i^n X')|)$$

$$\leq K A^3 \sqrt{\Delta_n} \sum_{i=1}^{[t/\Delta_n]} \sqrt{\mathbb{E}(\|\Delta_i^n X''\|^2 \wedge (A^2 \Delta_n))}$$

$$\leq K t A^3 \left( \frac{A\sqrt{\Delta_n}}{\theta} + \sqrt{\Gamma''(\theta)} \right)$$

for all $\theta \in (0, 1]$. Taking $\theta = \theta_n = \Delta_n^{1/4}$, so $\Gamma''(\theta_n) \to 0$, we deduce $\mathbb{E}(G_n) \to 0$; hence, we have (8.52).

*Step* 3. Observe that, for all $A > 1$ and $\rho > 0$, we have

$$\frac{1}{\Delta_n} V(f_{A\sqrt{\Delta_n}}, \Delta_n)_T \leq \widehat{A}'(\Delta_n)_T \leq \frac{1}{\Delta_n} V(f_\rho, \Delta_n)_T$$

as soon as $2A\sqrt{\Delta_n} \leq \alpha \Delta_n^{\varpi} \leq \rho$, that is for all $n$ large enough. We also, obviously, have (for each $\omega$)

$$\widetilde{D}(f_\rho)_T \xrightarrow{\rho \to 0} 0, \qquad \int_0^T \rho_{\sigma_s}(f_A) \, ds \xrightarrow{A \to \infty} C_T.$$

At this stage, (8.47) readily follows from (8.50) and (8.51). Hence, we are finished.

8.9. *Proof of Theorem 4.4.* Both (4.17) and (4.18) amount to proving

$$(8.53) \qquad V_n := \widetilde{\mathbb{P}}(|\widehat{H}_t^n| > Z_n \mid \mathcal{F}) \xrightarrow{\mathbb{P}} V := \widetilde{\mathbb{P}}(|H_t| > Z \mid \mathcal{F})$$

for $Z_n$ and $Z$ as in the statement of the theorem, and with

$$\widehat{H}_t^n = \sum_{i=1}^{[t/\Delta_n]} \sum_{j=1}^{2} g_n^j (\Delta_i^n X) \phi^j (\widehat{R}(n)_i), \qquad H_t = \sum_{q:\, S_q \le t} \sum_{j=1}^{2} g^j (\Delta X_{S_q}) \phi^j (\widehat{R}_q)$$

for the following choices of the $\mathbb{R}^2$-valued $g_n$, $\phi$, $\widehat{R}(n)_i$ and $\widehat{R}_q$: we refer to case 1 or 2, if we want to prove (4.18) or (4.17); then, we take $g_n(x) = g(x)\mathbf{1}_{\{\|x\| > \alpha \Delta_n^{\varpi}\}}$ and:

- for case 1: $g^j(x) = x_j^2, \phi^j(x) = x_{3-j}^2, \widehat{R}(n)_i = R(n)_i, \widehat{R}_q = R_q$,
- for case 2: $g^j(x) = x_j^2 x_{3-j}, \phi^j(x) = x_{3-j}, \widehat{R}(n)_i = R(n)_i', \widehat{R}_q = R_q'$.

*Step* 1. In a first step, we truncate the functions $g_n$ at some level $\rho > 0$, and the proof is somewhat similar to the proof of (4.11), whose notation is generally used, like $a_n(\rho)$ or $\Omega_n(t, \rho/2)$, and $S_q'$ and $i(n, q)$ when $\rho$ is fixed. We define a number of processes:

$$\widehat{H}(\rho)_t'^n = \sum_{i=1}^{[t/\Delta_n]} \sum_{j=1}^{2} (g_n^j \psi_\rho)(\Delta_i^n X) \phi^j (\widehat{R}(n)_i), \qquad \widehat{H}(\rho)^n = \widehat{H}^n - \widehat{H}(\rho)'^n,$$

$$H(\rho)_t^n = \sum_{q:\, S_q' \le t} \sum_{j=1}^{2} (g^j(1 - \psi_\rho))(\Delta X_{S_q'}) \phi^j (\widehat{R}(n)_{i(n,q)}),$$

$$H(\rho)'^n = \widehat{H}^n - H(\rho)^n.$$

$$H(\rho)_t = \sum_{q:\, S_q \le t} \sum_{j=1}^{2} (g^j(1 - \psi_\rho))(\Delta X_{S_q}) \phi^j (\widehat{R}_q), \qquad H(\rho)' = H - H(\rho).$$

We obviously have, for all $\rho, \varepsilon > 0$,

$$(8.54) \quad \begin{cases} \widetilde{\mathbb{P}}(|H(\rho)_t^n| > Z_n + \varepsilon \mid \mathcal{F}) - \widetilde{\mathbb{P}}(|H(\rho)_t'^n| > \varepsilon \mid \mathcal{F}) \\ \qquad \le V_n \\ \qquad \le \widetilde{\mathbb{P}}(|H(\rho)_t^n| > Z_n - \varepsilon \mid \mathcal{F}) + \widetilde{\mathbb{P}}(|H(\rho)_t'^n| > \varepsilon \mid \mathcal{F}), \\ \widetilde{\mathbb{P}}(|H(\rho)_t| > Z + \varepsilon \mid \mathcal{F}) - \widetilde{\mathbb{P}}(|H(\rho)_t'| > \varepsilon \mid \mathcal{F}) \\ \qquad \le V \\ \qquad \le \widetilde{\mathbb{P}}(|H(\rho)_t| > Z - \varepsilon \mid \mathcal{F}) + \widetilde{\mathbb{P}}(|H(\rho)_t'| > \varepsilon \mid \mathcal{F}). \end{cases}$$

*Step* 2. Observe that $|(g_n \psi_\rho)(x + y)| \le K(\|x\|^{m+2} \Delta_n^{-2\varpi} + (\|y\|^2 \wedge \rho^2))$ for $m = 2$ or $m = 3$, according to cases 1 or 2, and when $\rho \le 1/2$. Then, we deduce, exactly as for (8.37), that $\widetilde{\mathbb{E}}(|\widehat{H}(\rho)_t'^n|) \le a_n(\rho)$, and thus

$$(8.55) \qquad \lim_{\rho \to 0} \limsup_n \widetilde{\mathbb{E}}(|\widehat{H}(\rho)_t'^n|) = 0.$$

Next, we consider the set $B$ of all $\rho > 0$ such that outside a $\mathbb{P}$-null set we have $\|\Delta X_s(\omega)\| \notin B$ for all $s > 0$. This set $B$ has an at most countable complement. Suppose that $\rho$ is fixed in $B$. On the set $\Omega_n(t, \rho/2)$, we have

$$\widehat{H}(\rho)^n_t = \sum_{q \,:\, S'_q \leq \Delta_n[t/\Delta_n]} \sum_{j=1}^{2} (g^j_n(1-\psi_\rho))(\Delta^n_{i(n,q)}X)\phi^j(\widehat{R}(n)_{i(n,q)}).$$

We have $(g_n(1-\psi_\rho))(\Delta^n_{i(n,q)}X) \to (g(1-\psi_\rho))(\Delta X_{S'_q})$ a.s. for each $q$ (because $\rho \in B$). Recalling (8.22), we deduce that $\widehat{H}(\rho)^n_t - H(\rho)^n_t \xrightarrow{\widetilde{\mathbb{P}}} 0$ as $n \to \infty$. Combining this with (8.55), we obtain

$$(8.56) \qquad \lim_{\rho \to 0, \rho \in B} \limsup_n \widetilde{\mathbb{P}}(|H(\rho)'^n_t| > \eta) = 0 \qquad \forall \eta > 0.$$

*Step* 3. The variable $H(\rho)_t$ has the same form than $H(\rho)^n_t$, except that $\widehat{R}(n)_{i(n,q)}$ is substituted with $\widetilde{R}'_q$, the variable equal to $\widehat{R}_l$ on the set $\{S'_q = S_l\}$. There is no connection between $\widehat{R}(n)_{i(n,q)}$ and $\widetilde{R}'_q$. However, there is a strong connection between their $\mathcal{F}$-conditional laws, and hence between the $\mathcal{F}$-conditional laws $\zeta^n_\rho(\omega, dy)$ of $|H(\rho)^n_t|$ and $\zeta_\rho(\omega, dy)$ of $|H(\rho)_t|$.

More precisely, let $\sigma$ and $\sigma'$ be two $2 \times 2$ matrices with squares $c = \sigma\sigma^\star$ and $c' = \sigma'\sigma'^\star$. The variable $\phi(\sqrt{\kappa_1}\sigma U_1 + \sqrt{1 - \kappa_1}\sigma' U'_1)$ in case 1 or $\phi(\sqrt{L_1}\sigma U_1 + \sqrt{1 - k - L_1}\sigma' U'_1)$ in case 2 has a law $\theta_{c,c'}$ which depends only on $(c, c')$ and is continuous in $(c, c')$ (for the weak convergence). Moreover, since $U_1$ and $U'_1$ are independent standard two-dimensional Gaussian variables, the following property is immediate. If we have a sequence $(\Phi_q)$ of independent variables, each $\Phi_q$ being distributed according to $\theta_{c_q, c'_q}$, and if $(x(q))$ is a sequence of nonrandom vectors taking values in $\mathbb{R}^2$, then

$$(8.57) \qquad Q \geq 1, a > 0 \quad \Rightarrow \quad \mathbb{P}\left(\left|\sum_{q=1}^{Q}\sum_{j=1}^{2} x(q)_j \Phi^j_q\right| = a\right) = 0.$$

With this in mind, we see that $\zeta^n_\rho(\omega, .)$ is in fact the law of

$$\left| \sum_{q \,:\, S'_q(\omega) \leq t} \sum_{j=1}^{2} (g^j(1-\psi_\rho))(\Delta X_{S'_q}(\omega))\Phi^j_q \right|,$$

where the $\Phi_q$ are independent variables conditionally on $\mathcal{F}$, with the laws $\theta_{\widehat{c}(n,-)_{i(n,q)}(\omega), \widehat{c}(n,+)_{i(n,q)}(\omega)}$. Of course $\zeta_\rho$ is the same, with $\widehat{c}(n,-)_{i(n,q)}$ and $\widehat{c}(n,+)_{i(n,q)}$ substituted with $c_{S'_q-}$ and $c_{S'_q}$. Then, in view of (8.44) and of the continuity of $\theta_{c,c'}$ in $(c, c')$ and of (8.57), we readily deduce that, for any $\mathcal{F}$-measurable variables $Z_n \xrightarrow{\mathbb{P}} Z > 0$,

$$(8.58) \qquad \zeta^n_\rho((Z_n, \infty)) \xrightarrow{\mathbb{P}} \zeta_\rho((Z, \infty)).$$

*Step* 4. At this stage, we deduce from (8.54) that

$$
\begin{aligned}
|V_n - V| \leq\ & |\zeta_\rho^n((Z_n + \varepsilon, \infty)) - \zeta_\rho((Z + \varepsilon, \infty))| \\
& + |\zeta_\rho^n((Z_n - \varepsilon, \infty)) - \zeta_\rho((Z - \varepsilon, \infty))| \\
& + \zeta_\rho([Z - \varepsilon, Z + \varepsilon]) + \widetilde{\mathbb{P}}(|H(\rho)_t'^n| > \varepsilon \mid \mathcal{F}) + \widetilde{\mathbb{P}}(|H(\rho)_t'| > \varepsilon \mid \mathcal{F}).
\end{aligned}
$$

Then, in view of (8.58), and since $\{|H(\rho)_t - Z| \leq \varepsilon\} \subset \{|H_t - Z| \leq 2\varepsilon\} \cup \{|H(\rho)_t'| > \varepsilon\}$, we obtain

$$
\begin{aligned}
\limsup_n \mathbb{E}(|V_n - V|) \leq\ & \widetilde{\mathbb{P}}(|H_t - Z| \leq 2\varepsilon) \\
& + 2\widetilde{\mathbb{P}}(|H(\rho)_t'| \geq \varepsilon) + \limsup_n \widetilde{\mathbb{P}}(|H(\rho)_t'^n| \geq \varepsilon).
\end{aligned}
$$

This holds for all $\rho > 0$ and $\varepsilon > 0$. But the $\mathcal{F}$-conditional law of $|H_t|$ is either the Dirac mass $\varepsilon_0$ or it has a density, whereas $Z > 0$ by hypothesis, hence $\widetilde{\mathbb{P}}(|H_t - Z| \leq 2\varepsilon) \to 0$ as $\varepsilon \to 0$. Moreover $H(\rho)_t' \xrightarrow{\mathbb{P}} 0$ as $\rho \to 0$, hence by (8.56) we obtain $\limsup_n \mathbb{E}(|V_n - V|) = 0$, which implies (8.53).

8.10. *Proof of Theorem 4.6.* Let us prove (4.22). Set

$$
\Psi_n = \frac{1}{\Delta_n} \Phi_n^{(d)} \sqrt{V(g_1, \Delta_n)_T V(g_2, \Delta_n)_T}.
$$

By Theorems 4.1 and 4.2 and by (4.3), the variables $\Psi_n - \widehat{A}_n$ converge stably in law, in restriction to $\Omega_T^{(d)}$, to $\widetilde{D}_T$, for the two possible choices of $\widehat{A}_n$ whereas on $\Omega_T^{(d)}$ the $\mathcal{F}$-conditional law of $\widetilde{D}_T$ admits a density. Therefore,

$$
\mathbb{P}(A \cap \{\Psi_n > Z_n + \widehat{A}_n\}) \to \widetilde{\mathbb{P}}(A \cap \{\widetilde{D}_T > Z\}),
$$

which is (4.22). The proof of (4.21) is similar.

8.11. *Proof of Theorem 5.1.* Set $U_n = (\Phi_n^{(j)} - 1)/\widehat{V}_n^{(j)}$, and let $A \in \mathcal{F}$.

PROOF OF (b). We use (5.6), so we have

$$
\mathbb{P}(C_n^{(j)} \cap A) = \mathbb{P}(\{|U_n| \geq 1/\sqrt{\alpha}\} \cap A).
$$

Theorem 4.5(a) yields that $U_n$ converges stably in law, in restriction to $\Omega_T^{(j)}$, to a variable $U$ with $\mathbb{E}(U^2 \mid \mathcal{F}) = 1$; hence, if $A \subset \Omega_T^{(j)}$

$$
\limsup_n \mathbb{P}(C_n^{(j)} \cap A) \leq \widetilde{\mathbb{P}}(\{|U| \geq 1/\sqrt{\alpha}\} \cap A) \leq \alpha \mathbb{P}(A),
$$

where the last inequality follows from Bienaymé–Chebyshev applied to the conditional law of $U$ knowing $\mathcal{F}$. This clearly yields $\alpha^{(j)} \leq \alpha$. □

PROOF OF (a). Now, we use (5.3), so we have

$$\mathbb{P}(C_n^{(j)} \cap A) = \mathbb{P}(\{|U_n| \geq z_\alpha\} \cap A)$$

for each $A \in \mathcal{F}$. If $X$ and $\sigma$ do not jump together, Theorem 4.5(a) yields that $U_n$ converges stably in law, in restriction to $\Omega_T^{(j)}$, to a variable which is $\mathcal{N}(0,1)$ conditionally on $\mathcal{F}$. Then if $A \subset \Omega_T^{(j)}$ we have

$$\mathbb{P}(C_n^{(j)} \cap A) \to \widetilde{\mathbb{P}}(\{|U| \geq z_\alpha\} \cap A) = \alpha \mathbb{P}(A).$$

This yields (5.5), and hence (5.4) as well. □

PROOF OF (c). (1) By (5.7), we have

$$\mathbb{P}(C_n^{(j)} \cap A) = \mathbb{P}\left(A \cap \left\{\frac{|\Phi_n^{(j)} - 1| V(f, \Delta_n)_T}{\sqrt{\Delta_n}} > Z_n^{(j)}(\alpha)\right\}\right).$$

Hence, in view of (4.21), the property (5.5), and thus (5.4) as well, will follow if we prove that

(8.59) $\begin{cases} Z_n^{(j)}(\alpha) \xrightarrow{\mathbb{P}} Z(\alpha), & \text{in restriction to } \Omega_T^{(j)}, \\ \quad \text{where } Z(\alpha) \text{ is a positive and } \mathcal{F}\text{-measurable variable,} \\ \quad \text{with } \widetilde{\mathbb{P}}(|\widetilde{G}_T| > Z(\alpha) \mid \mathcal{F}) = \alpha \text{ in restriction to } \Omega_T^{(j)}. \end{cases}$

By (3.11) and (3.15) and also (e) of Assumption (H), we see that the law of $|\widetilde{G}_T|$, conditional on $\mathcal{F}$ and in restriction to $\Omega_T^{(j)}$, has no atom (and indeed has a positive density on $\mathbb{R}_+$), so $Z(\alpha)$ satisfying $\widetilde{\mathbb{P}}(|\widetilde{G}_T| > Z(\alpha) \mid \mathcal{F}) = \alpha$ is uniquely defined and positive a.s. on $\Omega_T^{(j)}$. By (4.16), it is also obvious that the law of $|\widehat{G}_T^n|$, conditional on $\mathcal{F}$, has no atom except possibly $\{0\}$: hence if $0 < \gamma < 1$, the variable $Z_n'(\gamma) = \sup(a : \widetilde{\mathbb{P}}(|\widehat{G}_T^n| > a \mid \mathcal{F}) \geq \gamma)$ satisfies $\widetilde{\mathbb{P}}(|\widehat{G}_T^n| > Z_n'(\gamma) \mid \mathcal{F}) \leq \gamma$, with equality as soon as $Z_n'(\gamma) > 0$. Then, (4.17) applied with $Z_n = Z = \gamma$, yields

(8.60) $\quad \gamma \in (0,1) \quad \Longrightarrow \quad Z_n'(\gamma) \xrightarrow{\mathbb{P}} Z(\gamma) \qquad \text{on the set } \Omega_T^{(j)}.$

(2) Consider an i.i.d. sequence of positive variables $Y_i$ with a purely nonatomic law, and denote by $Z$ the unique (decreasing) function such that $\mathbb{P}(Y_i > Z(x)) = x$ for all $x \in (0,1)$. We set $U_n(x) = \frac{1}{N_n} \sum_{i=1}^{N_n} 1_{\{Y_i > x\}}$, and call $V_n(\alpha)$ the $[\alpha N_n]$th variable, after they have been rearranged in decreasing order, for some $\alpha \in (0,1)$. Assume that $N_n > 4/\alpha(1-\alpha)$ and take $\varepsilon \in (4/N_n, \alpha(1-\alpha))$. If $V_n(\alpha) > Z(\alpha - \varepsilon)$, we have

$$U_n(Z(\alpha - \varepsilon)) \geq U_n(V_n(\alpha)) = \frac{[\alpha N_n] - 1}{N_n} \geq \alpha - \frac{2}{N_n} \geq \alpha - \frac{\varepsilon}{2},$$

that is, $U_n(Z(\alpha - \varepsilon)) - (\alpha - \varepsilon) \geq \varepsilon/2$. In a similar way, if $V_n(\alpha) < Z(\alpha + \varepsilon)$ we have $U_n(Z(\alpha + \varepsilon)) - (\alpha + \varepsilon) \leq -\varepsilon$. Since the variables $U_n(Z(x))$ have mean $x$ and variance smaller than $1/4N_n$, it follows from the Bienaymé–Chebyshev inequality that

$$
\begin{aligned}
(8.61) \qquad & \mathbb{P}(V_n(\alpha) \notin [Z(\alpha + \varepsilon), Z(\alpha - \varepsilon)]) \\
& \leq \mathbb{P}(U_n(Z(\alpha - \varepsilon)) - (\alpha - \varepsilon) \geq \varepsilon/2) \\
& \quad + \mathbb{P}(U_n(Z(\alpha + \varepsilon)) - (\alpha + \varepsilon) \leq \varepsilon) \\
& \leq \frac{5}{4N_n \varepsilon^2}.
\end{aligned}
$$

(3) Now, recall from (5.2) that $Z_n^{(j)}(\alpha)$ is the $[\alpha N_n]$th absolute order statistics for $N_n$ independent draws of $\widehat{G}(\Delta_n)_T$, conditionally on $\mathcal{F}$. Then, (8.61) with the choice $\varepsilon = 1/N_n^{1/4}$, together with (8.60), imply that

$$
\widetilde{\mathbb{P}}(Z_n^{(j)}(\alpha) \notin [Z(\alpha - N_n^{-1/4}) - \eta, Z(\alpha + N_n^{-1/4}) + \eta] \mid \mathcal{F}) \xrightarrow{\mathbb{P}} 0
$$

for all $\eta > 0$, on the set $\Omega_T^{(j)}$, and (8.59) follows. $\square$

8.12. *Proof of Theorem 5.2.* (1) The first thing is that the new cutoffs give us the same level as the old ones. But, on the set $\Omega_T^{(j)}$, we know that $\widehat{V}_n^{(j)}/\sqrt{\Delta_n}$ converges in probability to a finite limit; hence,

$$
\mathbb{P}(\{\widehat{V}_n^{(j)} = (\widehat{V}_n^{(j)} \wedge (\alpha' \Delta_n^{\varpi'}))\} \cap \Omega_T^{(j)}) \to \mathbb{P}(\Omega_T^{(j)})
$$

and the result is obvious.

(2) The second thing we have to prove is $\beta^{(j)} = 1$. This amounts to proving that $\mathbb{P}(C_n^{(j)} \cap \Omega_T^{(d)}) \to \mathbb{P}(\Omega_T^{(d)})$, or equivalently that $\mathbb{P}(\{|U_n| \geq \eta\} \cap \Omega_T^{(d)}) \to \mathbb{P}(\Omega_T^{(d)})$ for any fixed $\eta > 0$, where here $U_n = (\Phi_n^{(j)} - 1)/(\widehat{V}_n^{(j)} \wedge (\alpha' \Delta_n^{\varpi'}))$. Since $\Phi_n^{(j)} - 1$ converges stably in law to an almost surely nonvanishing limit by Theorem 3.1(b), on $\Omega_T^{(d)}$, the result will be implied by the property that $(\widehat{V}_n^{(j)} \wedge (\alpha' \Delta_n^{\varpi'})) \to 0$, which is obvious.

8.13. *Proof of Theorem 5.3.* The proof is the same as for Theorem 5.1, with a few changes. In case (a), that is of (5.11), we set $U_n = \Phi_n^{(d)}/\widehat{V}_n$, so

$$
\mathbb{P}(C_n^{(d)} \cap A) = \mathbb{P}(\{|U_n| \geq 1/\alpha\} \cap A).
$$

Theorem 4.5(b) yields that $U_n$ converges stably in law, in restriction to $\Omega_T^{(d)}$, to a variable $U > 0$ with $\widetilde{\mathbb{E}}(U \mid \mathcal{F}) = 1$; hence, if $A \subset \Omega_T^{(d)}$,

$$
\limsup_n \mathbb{P}(C_n^{(d)} \cap A) \leq \widetilde{\mathbb{P}}(\{|U| \geq 1/\alpha\} \cap A) \leq \alpha \mathbb{P}(A),
$$

(use again the Markov inequality), and thus $\alpha^{(d)} \leq \alpha$. The property $\beta^{(d)} = 1$ amounts to having $\mathbb{P}(\{U_n \geq \eta\} \cap \Omega_T^{(j)}) \to \mathbb{P}(\Omega_T^{(j)})$ for any fixed $\eta > 0$. By Theorem 3.1(a), we have $\Phi_n^{(d)} \xrightarrow{\mathbb{P}} B_T/\sqrt{B_T'^1 B_T'^2} > 0$ on $\Omega_T^{(j)}$, and also $\widehat{V}_n \xrightarrow{\mathbb{P}}$ 0 on this set by (4.3), so $U_n \xrightarrow{\mathbb{P}} +\infty$ on $\Omega_T^{(j)}$, and the result readily follows.

Finally, in case (b) of (5.13), the proof is exactly the same as for case (c) of Theorem 5.1.

**Acknowledgments.** We would like to thank an associate editor and an anonymous referee for very constructive comments which greatly improved the paper.

Institut de Mathématiques de Jussieu
CNRS–UMR 7586
Université Pierre et Marie Curie–P6
4 Place Jussieu, 75252 Paris-Cedex 05
France
E-mail: jean.jacod@upmc.fr

Department of Finance
Kellogg School of Management
Northwestern University
Evanston, Illinois 60208-2001
USA
E-mail: v-todorov@kellogg.northwestern.edu